\documentclass{article}
\usepackage[utf8]{inputenc}
\usepackage{titlesec}

\titleformat{\section}
{\normalfont\normalsize\bfseries}{\thesection}{1ex}{}
\titleformat{\subsection}
{\normalfont\normalsize\bfseries}{\thesubsection}{1ex}{}
\titleformat{\subsubsection}[runin]
{\normalfont\normalsize\bfseries}{\thesubsubsection}{1ex}{}
\titleformat{\paragraph}[runin]
{\normalfont\normalsize\bfseries}{\theparagraph}{1ex}{}
\titleformat{\subparagraph}[runin]
{\normalfont\normalsize\bfseries}{\thesubparagraph}{1ex}{}

\titlespacing{\section}{0pt}{12pt plus 4pt minus 2pt}{0pt plus 2pt minus 2pt}
\titlespacing{\subsection}{0pt}{12pt plus 4pt minus 2pt}{0pt plus 2pt minus 2pt}
\titlespacing{\subsubsection}{0pt}{12pt plus 4pt minus 2pt}{0pt plus 2pt minus 2pt}

\usepackage{amsmath}
\usepackage{amssymb}
\usepackage{amsthm}

\usepackage{titling}
\thanksmarkseries{arabic}
\usepackage{caption}
\usepackage{subcaption}

\usepackage{geometry}
\usepackage{hyperref}
\usepackage{siunitx}

\usepackage{graphicx}
\usepackage{tikz}
\usepackage{placeins}
\usepackage{enumitem}

\usepackage{algorithm}
\usepackage[noend]{algpseudocode}

\setlist[itemize]{noitemsep, topsep=.5ex}
\setlist[enumerate]{noitemsep, topsep=.5ex}

\usepackage[backend = biber,
    style=numeric,
    giveninits = true,
    natbib = true,
    hyperref = true,
    maxbibnames = 4,
    sorting=none
    ]{biblatex}
    
\newlength{\MyWidth}
\newlength{\MyUnit}
\newcommand{\tr}{ {\mathsf T} }
\newcommand{\ct}{ {\mathsf H} }
\newcommand{\ima}{\operatorname{i}}
\newcommand{\mean}{\operatorname{mean}}
\newcommand{\diag}{\operatorname{diag}}
\newcommand{\env}{\operatorname{env}}
\newcommand{\vek}{\operatorname{vec}}
\newcommand{\nrm}{\operatorname{nrm}}

\newcommand{\DLT}{\operatorname{DLT}_\zeta} 
\newcommand{\FFT}{\operatorname{FFT}} 

\newcommand{\Sim}{\operatorname{sim}}
\newcommand{\Meas}{\operatorname{meas}}

\newcommand{\vsigma}{\pmb{\sigma}}
\newcommand{\vepsi}{\pmb{\varepsilon}}
\newcommand{\vbtau}{\pmb{\bar{\tau}}}
\newcommand{\vgamma}{\pmb{\gamma}}
\newcommand{\vomega}{\pmb{\omega}}

\newcommand{\vpsi}{\hat{\mathbf{\psi}}}
\newcommand{\mPsi}{\hat{\mathbf{\Psi}}}

\newcommand{\vnabla}{\pmb{\nabla}}

\newcommand{\mA}{\mathbf{A}}
\newcommand{\mC}{\mathbf{C}}
\newcommand{\mE}{\mathbf{E}}
\newcommand{\mF}{\mathbf{F}}
\newcommand{\mM}{\mathbf{M}}
\newcommand{\mI}{\mathbf{I}}
\newcommand{\mS}{\mathbf{S}}
\newcommand{\mU}{\mathbf{U}}
\newcommand{\mV}{\mathbf{V}}
\newcommand{\mX}{\mathbf{X}}

\newcommand{\fbu}{\boldsymbol{\bar{u}}}
\newcommand{\fbv}{\boldsymbol{\bar{v}}}
\newcommand{\fn}{\boldsymbol{n}}

\newcommand{\fub}{\hat{\boldsymbol{u}}_b}
\newcommand{\fN}{\boldsymbol{N}}

\newcommand{\va}{\mathbf{a}}
\newcommand{\vc}{\mathbf{c}}
\newcommand{\vh}{\mathbf{h}}
\newcommand{\vk}{\mathbf{k}}
\newcommand{\vo}{\mathbf{0}}

\newcommand{\vv}{\mathbf{v}}
\newcommand{\vy}{\mathbf{y}}
\newcommand{\vx}{\mathbf{x}}

\newcommand{\vs}{\mathbf{s}}
\newcommand{\vt}{\mathbf{t}}

\newcommand{\CC}{\mathbb{C}}
\newcommand{\RR}{\mathbb{R}}

\newenvironment{strech}{}{}

\newcommand{\EE}{\operatorname{E}}   
\newcommand{\pp}{\operatorname{p}}   
\newcommand{\NDoF}{N}                
\newcommand{\NMeas}{K}               
\newcommand{\NEval}{M}               
\newcommand{\Fo}{\mathbf{F}}         
\newcommand{\para}{q}                
\newcommand{\vpara}{\mathbf{q}}      
\newcommand{\paraspace}{\mathcal{Q}} 
\newcommand{\wn}{\hat k}             
\newcommand{\wl}{\hat \lambda}       
\newcommand{\nes}{s}                 
\newcommand{\ecf}{f_c}               
\newcommand{\ect}{t_c}               
\newcommand{\dOmega}{\omega_\Delta}  
\newcommand{\EWM}{\zeta}             
\newcommand{\SC}{\vx_c}              
\newcommand{\xo}{\chi}               
\newcommand{\NN}{\mathbb{N}}    
\title{Experimental validation of an inverse method for defect reconstruction in a 2D waveguide model}

\date{\vspace{-10ex}}

\author{J. Bulling\thanks{Bundesanstalt für Materialforschung und -pruefung, Unter den Eichen 87, 12205 Berlin, Germany \newline Jannis.Bullling@BAM.de, Jens.Prager@BAM.de} \and B. Jurgelucks\thanks{Department of Mathematics, Humboldt-Universität zu Berlin, Unter den Linden 6, 10099 Berlin, Germany\newline Benjamin.Jurgelucks@math.hu-berlin.de, Andrea.Walther@math.hu-berlin.de} \and J. Prager$^1$ \and A. Walther$^2$}

\begin{document}

\maketitle

\paragraph*{Abstract}

Defect reconstruction is essential in non-destructive testing and structural health monitoring with guided ultrasonic waves. This paper presents an algorithm for reconstructing notches in steel plates which can be seen as artificial defects representing cracks  by comparing measured results with those from a simulation model. The model contains a parameterized notch, and its geometrical parameters are to be reconstructed. While the algorithm is formulated and presented in a generalized form for many different  defect types, a special case of guided wave propagation is used to investigate one of the simplest possible simulation models that discretizes only the cross-section of the steel plate. An efficient simulation model of the plate cross-section is obtained by the semi-analytical Scaled Boundary Finite Element Method.  The reconstruction algorithm applied is gradient-based, and Algorithmic Differentiation calculates the gradient.  The dedicated experimental setup excites nearly plane wave fronts propagating orthogonal to the notch. A scanning Laser Doppler Vibrometer records the velocity field at certain points on the plate surface as input to the reconstruction algorithm. Using two plates with notches of different depths, it is demonstrated that accurate geometry reconstruction is possible.


\section{Introduction}
Characterizing damage in shell-like structures, such as plates, laminates, or pipes, is an important issue in non-destructive testing and structural mealth monitoring. 
Methods based on ultrasonic guided waves are appealing for damage characterization because the guided wave propagation allows the inspection of large areas of the structure from only a few sensor positions.
However, determining the damage size using these methods can be challenging, as the various dispersive modes of guided waves do not always permit a straightforward correlation between reflections from damage and damage size~\cite{lowe2002lowA0,lowe2002lowS0,demma2003reflection,demma2004reflection}.
In addition, the wave packet reflected from the defect must be distinguished from the natural wave propagation in the pristine structure for many methods. Therefore, many state-of-the-art methods are based on residual signals that compare a measurement of the current, possibly damaged, condition with a baseline measurement of the pristine condition~\cite{ihn2008pitch,zhao2011ultrasonic,loendersloot2015damage}. In order to use these methods, the baseline measurement of the pristine condition must be available, and it must be guaranteed that the pristine condition is free of damage. This requirement can be restrictive, and baseline-free alternatives are desirable.

A possible alternative and baseline-free approach to solve this problem is an algorithm that compares the measured signals with those generated by a simulation model with damage of a specific size. If the simulated signal is close to the measured signal, the damage in the simulation model is a reconstruction of the actual damage.

The main objective of this work is to validate an algorithm whose basic concept has been published in~\cite{bulling2022defect}. The algorithm is able to reconstruct a defect defined by certain geometric parameters.
The proposed algorithm is derived by solving an inverse problem formulated as an optimization problem. A forward model is iteratively updated until the output matches the measured signals as closely as possible. 

The first key aspect of our work is the application of an efficient forward model. The forward model is based on the semi-analytical Scaled Boundary Finite Element Method (SBFEM), first introduced by Song and Wolf~\cite{song1997scaled, wolf2000scaled, song2000scaled}. The SBFEM approximates the solution of the wave equations by dividing the domain into smaller subdomains called super-elements.
The solution approach for each super-element is derived by a coordinate transformation that uses a finite element approximation of parts of the boundary. 
Since this finite element boundary approximation contains the degrees of freedom of the approximation, one obtains an order of magnitude reduction in the degrees of freedom compared to a pure finite element approximation of the entire domain.
In this work, two types of super-elements are considered: On the one hand, the high-order approach for polygonal elements~\cite{bazyar2008continued,gravenkamp2012numerical,chen2014high} and on the other hand, the waveguide super-elements~\cite{gravenkamp2015simulation,gravenkamp2018efficient}.

The second key aspect is the gradient-based optimization procedure. The gradient is computed by applying Algorithmic Differentiation (AD)~\cite{GrWa08,Nau12} to the forward model. By computing the gradient using AD, the optimization is fast and accurate. In addition, the AD framework allows for a flexible algorithm structure that can be used to model different objective functions, boundary conditions, and damage types~\cite{bulling2022defect}.

The work presented deals with the reconstruction of rectangular notches in steel plates using planar ultrasonic guided waves.
Investigations of notches as a surrogate for natural cracks in guided wave studies have a long tradition~\cite{lowe2002lowA0,lowe2002lowS0,demma2003reflection,demma2004reflection,wang2010evaluation,brence2022determination}.
While the description of the proposed algorithm can be easily extended to general problems, this artificial defect is used because it can be manufactured with high precision.
The wave propagation and the type of notches studied in the current work can also be mathematically reduced to a two-dimensional problem by considering only the cross-section of the plate.
Reducing the mathematical model to the cross-section often allows much faster computation, and the problem is more limited in the number of parameters. On the other hand, the experiment must be designed to fit the mathematical model of propagating waves. This paper, therefore, discusses the challenges of validating 2D models through experimentation. It should be emphasized that the validation of 2D models is a crucial step for the evaluation of a general defect reconstruction method.

Regarding model-based reconstruction, an alternative to geometric optimization of defects, defects can be modeled as material changes. This approach is known as Full Waveform Inversion (FWI), and was pioneered for seismic waves~\cite{fichtner2010full}. The references~\cite{rao2016guided, rao2017investigation, rao2020multi} are examples of FWI for guided waves that include experimental validation. A similar idea can be found in~\cite{huthwaite2013high}, where wave velocities define the model parameters.
Although these results are very promising, the methods require a numerical model that allows for a material change at each point where a defect might exist, making these models computationally expensive. In contrast, the method presented here takes advantage of the fact that large parts of the structure can be calculated semi-analytically by assuming a homogeneous material model, resulting in more efficient models.

This paper is organized as follows. Section~\ref{sec:Data} summarizes the measurement setup for the two investigated plates. Section~\ref{sec:Theory} defines the proposed defect reconstruction algorithm. Section~\ref{sec:Results} analyzes the performance of the algorithm and shows the reconstruction results.
The conclusion is given in section~\ref{sec:Conclusion}.


\section{Experimental Set-up and Data Acquisition}~\label{sec:Data}

\begin{figure}[htb]
    \centering
    \subcaptionbox{Photo of the two plates with notch depth \SI{0.4}{\milli\meter} and \SI{0.8}{\milli\meter}.}
    {\includegraphics[width=\MyWidth]{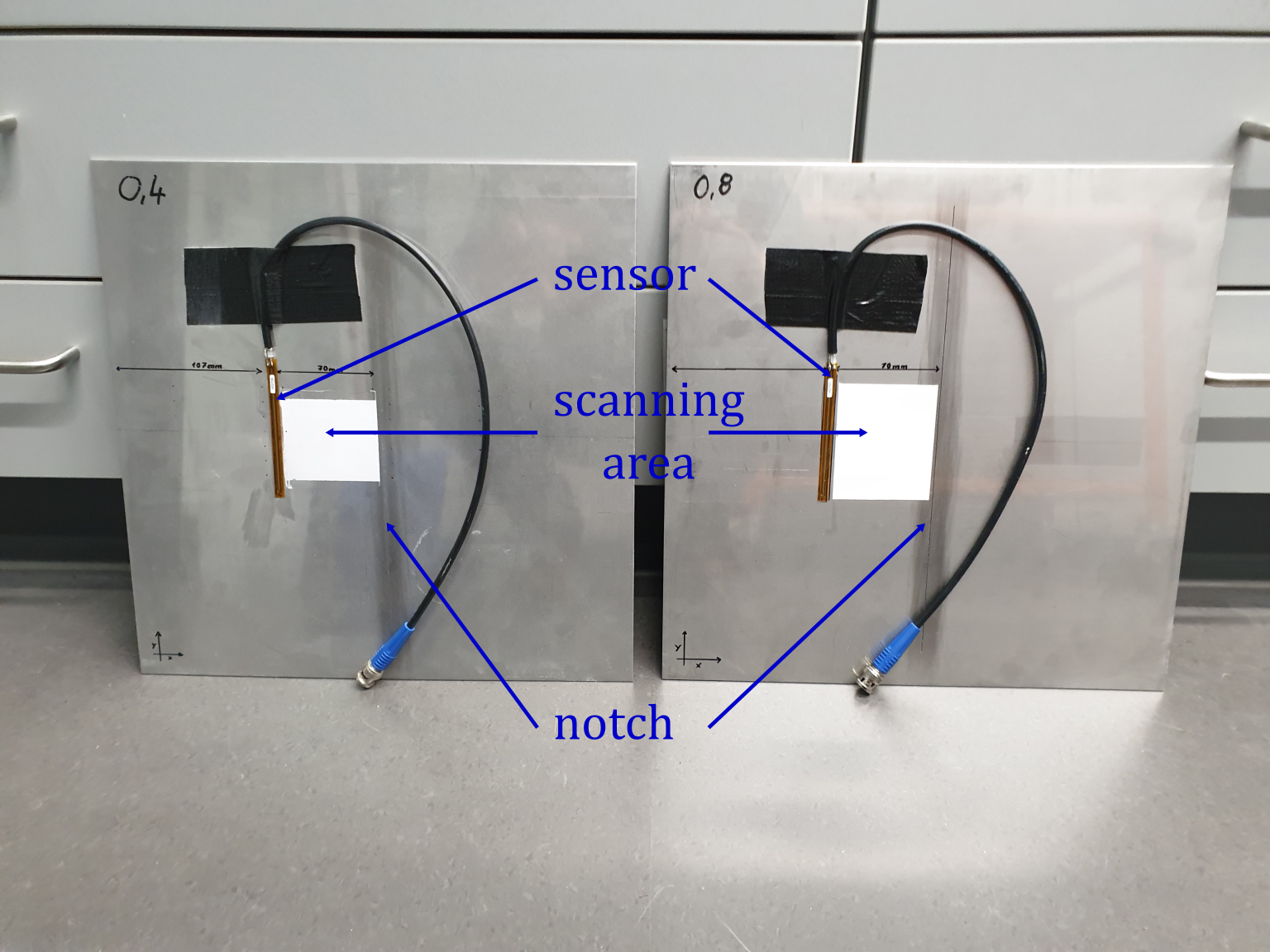}}
    \hfil
    \subcaptionbox{Photo of the experimental set-up.}
    {\includegraphics[width=\MyWidth]{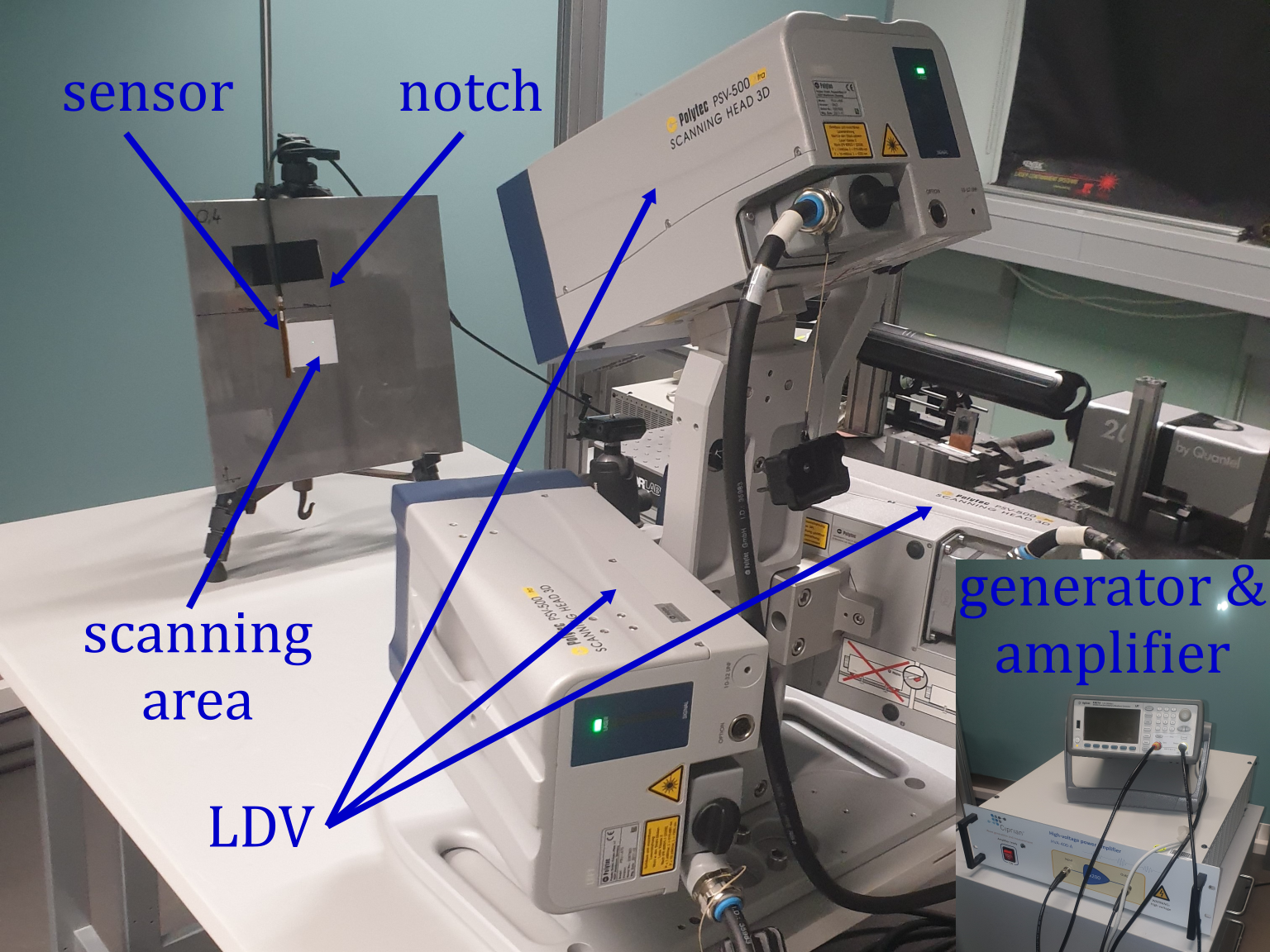}}
    \caption{Photos in the laboratory.}
    \label{fig:Photos}
\end{figure}

\begin{figure}[ht!]
    \centering
    \begin{tikzpicture}[x=\MyUnit,y=\MyUnit]
    \node[inner sep=0pt] (Overview) at (-4.5,0)
        {\includegraphics[width=5\MyUnit]{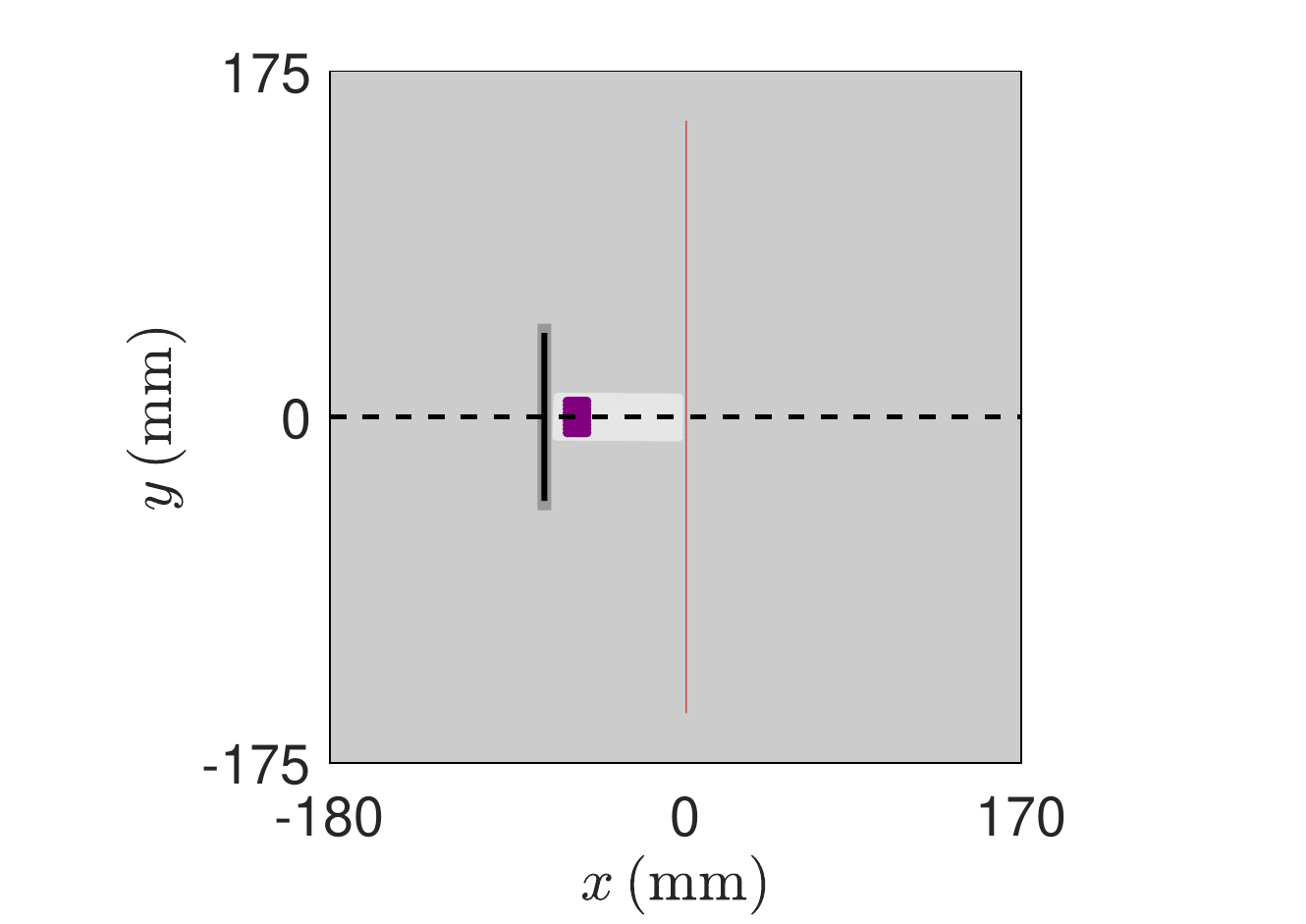}};
    \node[inner sep=0pt] (Detail1) at (-1.25,.5)
        {\fbox{\includegraphics[width=3\MyUnit]{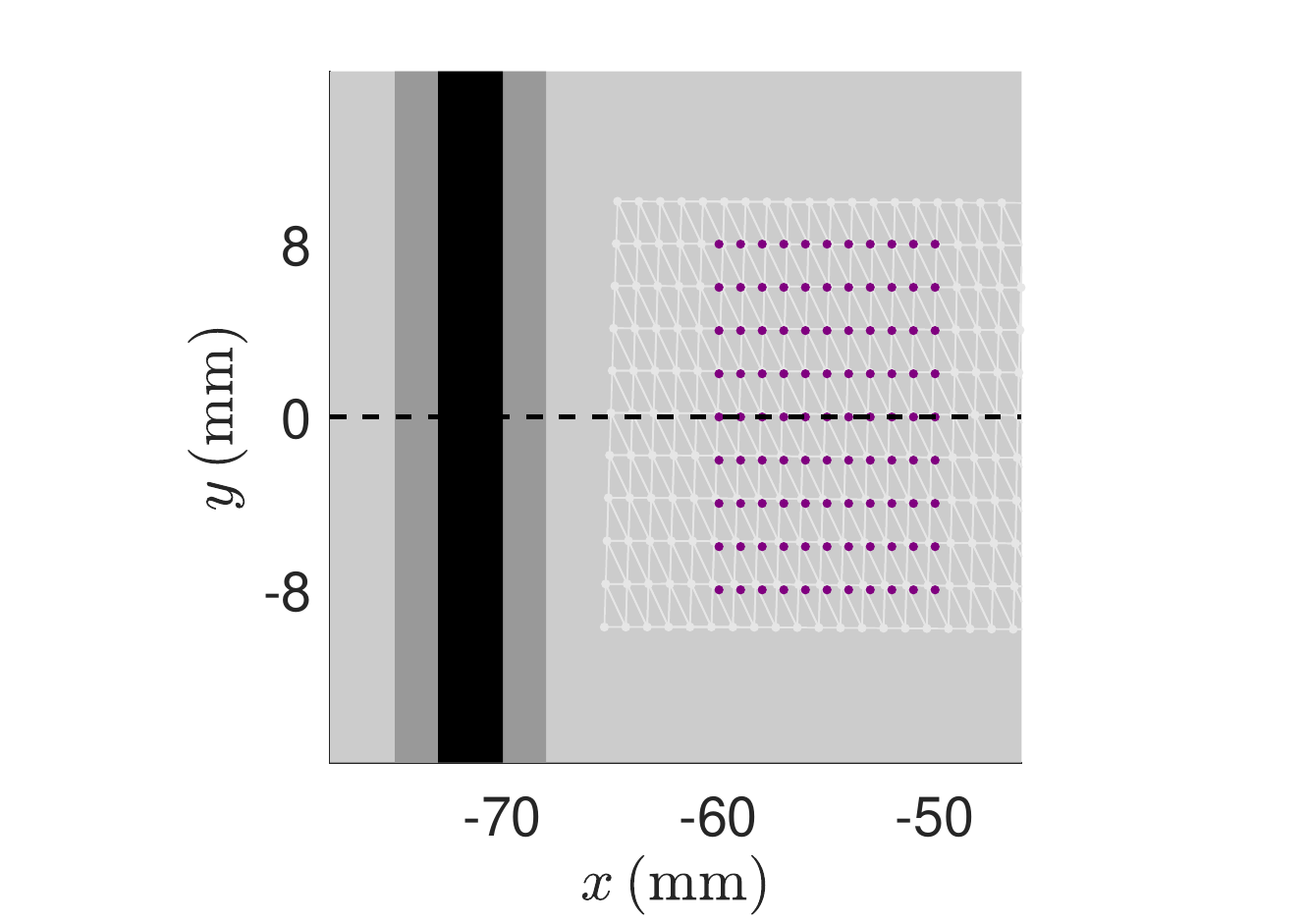}}};
    \draw[-,thick] (-4.8,.2) -- (Detail1.west);
    \node[inner sep=0pt] (Detail2) at (1.25,-.5)
        {\includegraphics[width=2.5\MyUnit]{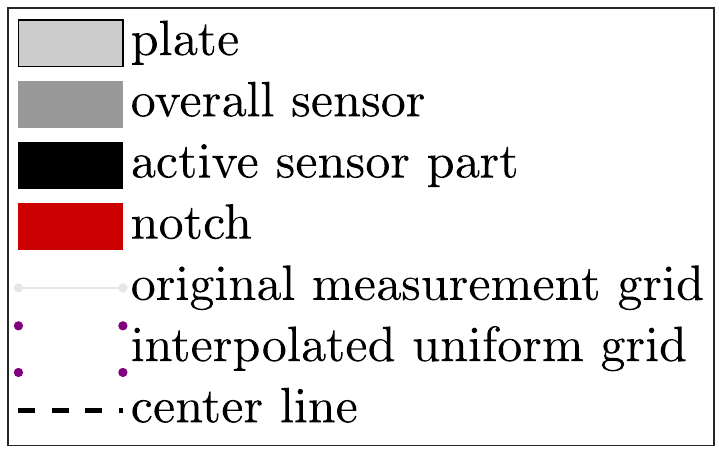}};
    \end{tikzpicture}
    \caption{Overview of the LDV measurement of the plate with the original measurement grid of the case with the notch depth $\SI{0.8}{\milli\meter}$.\label{fig:LDVoverview}}
\end{figure}

The reconstruction is studied on two $\SI{350}{\milli\meter} \times \SI{350}{\milli\meter} \times \SI{2}{\milli\meter}$ steel plates having a notch near the center.
The notches are $\SI{0.5}{\milli\meter}$ wide and $\SI{300}{\milli\meter}$ long, while the notch depth differs with $\SI{0.4}{\milli\meter}$ and $\SI{0.8}{\milli\meter}$, respectively.
Figure~\ref{fig:Photos} shows photos of the plates and an overview of the experimental set-up, while Figure~\ref{fig:LDVoverview} presents the schematic layout of the plates.
The positioning of the origin of the coordinate system is arbitrary, but for convenience, the origin was chosen to be at the left edge of the notch. This makes deviations of the reconstruction from the nominal notch position trivial to calculate.

As mentioned in the introduction, the reconstruction algorithm is based on a 2D cross-sectional model. However, a cross-sectional model simulates only plane wave fronts, which are constant orthogonal to the direction of wave propagation. Therefore, the experimental setup should favor the propagation of plane waves to achieve the best possible agreement between the simulation model and the measurement.

In summary, the data collection process can be divided into the following steps, which will be discussed in more detail below:
\begin{enumerate}
    \item The excitation signal is defined on a PC.
    \item Based on step 1, the electrical input signal for the sensor is generated and amplified.
    \item The sensor uses the electrical signal to excite the ultrasonic waves inside the plate.
    \item A laser Doppler vibrometer (LDV) measures the ultrasonic waves across a scanning area (Figure~\ref{fig:Photos}).
    \item The measurement is post-processed to create the input for the reconstruction algorithm.
\end{enumerate}

\begin{figure}[ht!]
    \centering
    \subcaptionbox{Normalized excitation signal and its envelope in the time-domain.\label{fig:excitationA}}
    {\includegraphics[width=\MyWidth]{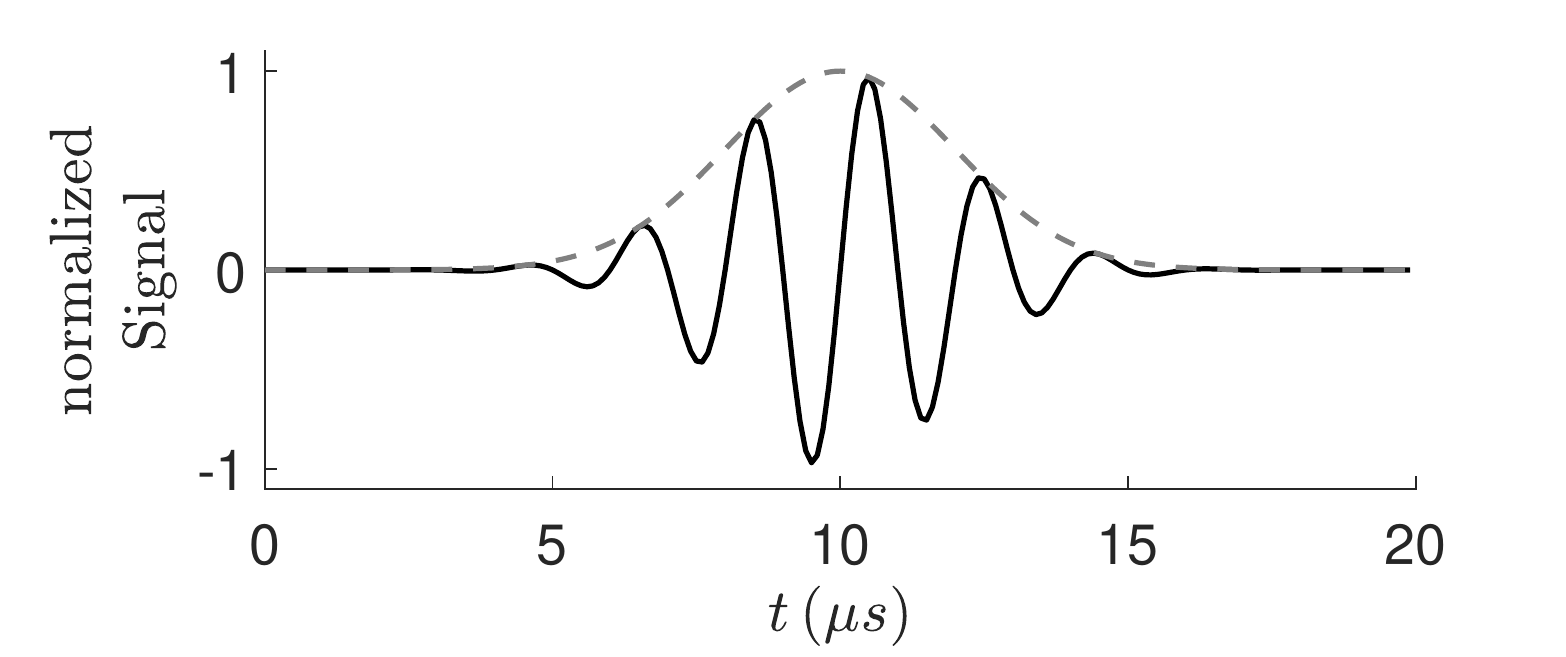}}
    \hfil
    \subcaptionbox{Wavelength of the modes with A1 cutoff frequency (dashed line) and normalized excitation signal in the frequency-domain.\label{fig:excitationB}}
    {\includegraphics[width=\MyWidth]{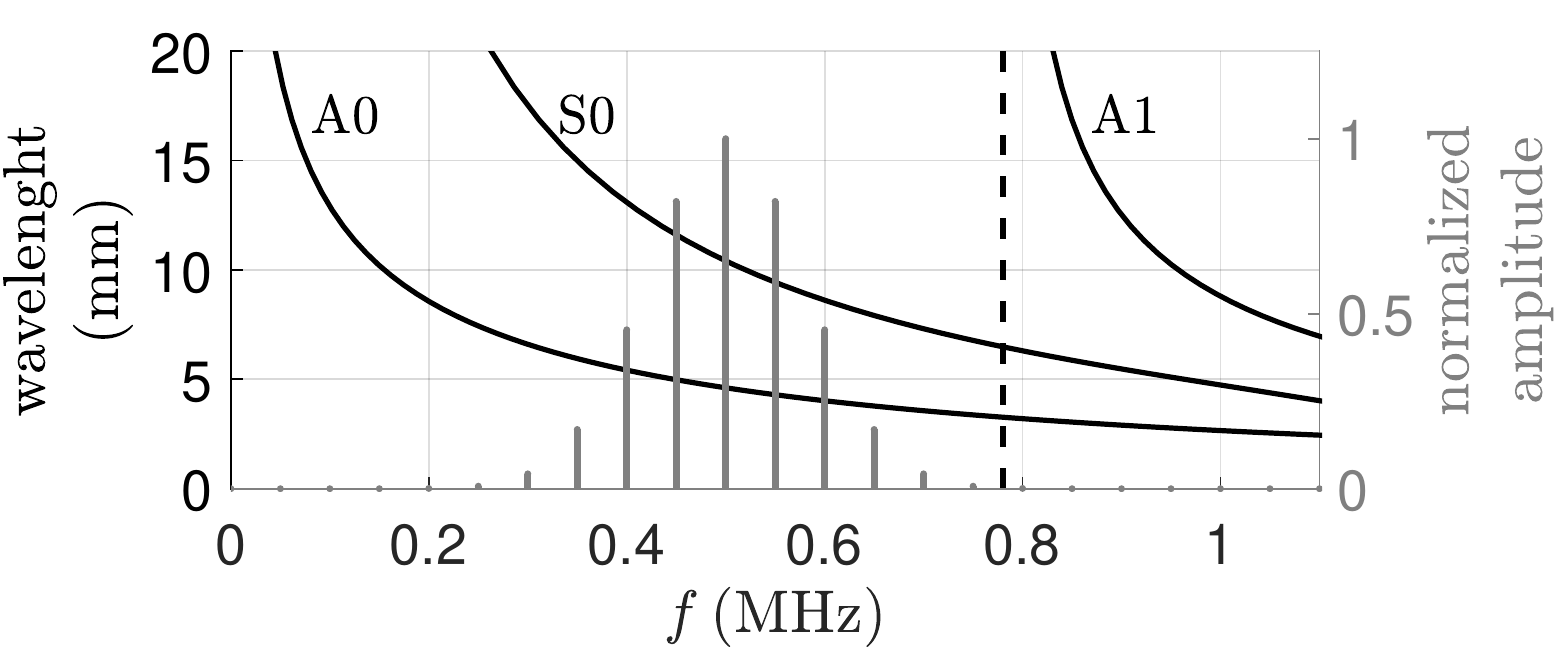}}
    \caption{Time and frequency dependency of the excitation signal.}
    \label{fig:excitation}
\end{figure}

\textbf{Step 1:}  Data acquisition begins with the definition of the normalized excitation signal $\nes$. 
Here, the normalized excitation signal is a sinusoidal Gaussian pulse, i.e.
$\nes(t) = \sin(2\pi \ecf t)\exp\left(- 0.5{(t - \ect)^2}/{\ecf^{-2}}\right)$
with a center frequency $\ecf$ of $\SI{500}{\kilo\hertz}$ and a time shift $\ect$ of $5 \ecf^{-1} = \SI{10}{\micro\second}$.
The choice of a Gaussian pulse is intended to promote agreement between experiment and simulation, since numerical models for smooth excitations achieve higher accuracy with the same number of degrees of freedom. Figure~\ref{fig:excitationA} shows the normalized excitation signal in the time domain, while Figure~\ref{fig:excitationB} shows the amplitude in the frequency domain with renormalization. In addition, Figure~\ref{fig:excitationB} shows the wavelengths of the S0, A0 and S1 modes of the plate. For the corresponding spectrum, only the S0 and A0 modes are excitable. The dispersion relation between frequency and wavelength is calculated by the SBFEM~\cite{gravenkamp2012numerical,bulling2020sensitivity}.

\textbf{Step 2 and 3:} The normalized excitation signal $\nes$ is converted into an electrical signal using the 
arbitrary waveform generator
$\langle$33521A, Agilent/dataTec GmbH, Reutlingen, Germany$\rangle$
and then amplified by a high-voltage amplifier 
$\langle$HVA-400-A, Ciprian, Grenoble, France$\rangle$.
The amplified signal is sent to a Macro Fiber Composite sensor (MFC-sensor) 
$\langle$M8503P2MFC, Smart Material GmbH, Dresden, Germany$\rangle$.
The MFC-sensor has an active piezoelectric area of $\SI{3}{\milli\meter} \times \SI{85}{\milli\meter}$, and a total area of $\SI{7}{\milli\meter} \times \SI{94.5}{\milli\meter}$.  
The long extension of the sensor in the $y$-direction compared to the wavelength of the propagating modes (about $\wl_{A0} = \SI{4.6}{\milli\meter}$ and $\wl_{S0} = \SI{10.4}{\milli\meter}$) produces a nearly planar wavefront and thus to the wave type that also fits the model assumption of a plane strain state. 
The velocities in the $y$-direction are decoupled from the velocities in the $x$-and $z$-directions for plane waves propagating in the $x$-direction. Due to the choice of the sensor, SH modes associated with oscillations in the $y$-direction are negligibly excited, so only velocities in the $x$-and $z$-directions are considered.

\begin{figure}[t]
    \centering
    \subcaptionbox{$v_x$ (in-plane)}
    {\includegraphics[height=\MyWidth, angle=-90]{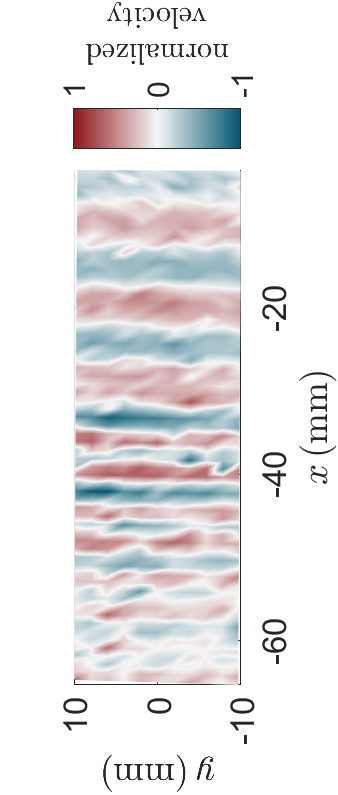}}
    \hfil
    \subcaptionbox{$v_z$ (out-of-plane)}
    {\includegraphics[height=\MyWidth, angle=-90]{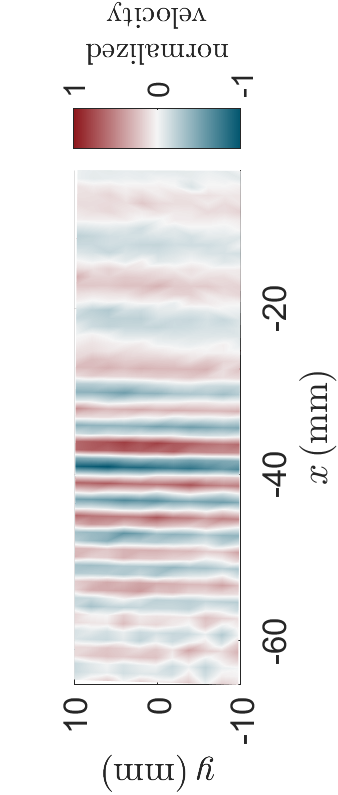}}
    \caption{Normalized velocity fields at $t=\SI{20}{\micro\second}$ for the original measurement grid (Figure~\ref{fig:LDVoverview}) on the plate with a notch of $\SI{0.8}{\milli\meter}$ depth.}
    \label{fig:XYwavefield}
\end{figure}

\textbf{Step 4:} A scanning LDV
$\langle$PSV-500-3D-HV, Polytec GmbH, Waldbronn, Germany$\rangle$
records the resulting ultrasonic waves. Due to the high reflectivity of the steel plate, it was necessary to treat the scanning area with white, non-permanent paint to ensure a defused reflection of the laser (see Figure~\ref{fig:Photos}). 
The scanning LDV measures the velocity field in all three directions in a predefined grid.
For our application, the $y$-velocity field is neglected.
To increase the signal-to-noise ratio, the signal is averaged over $\num{1000}$ times for each measurement point.
As mentioned above, the model simulates planar waves. Therefore, an evaluation of the planar waves from the experimental data is necessary for a good fit. Figure~\ref{fig:XYwavefield} shows that the excitation already produces a nearly planar wave. Nevertheless, some small changes in the $y$-direction can be observed. As a filter for waves with a wave vector not parallel to the $x$-axis, the average of the velocities along the $y$-axis is taken.

\textbf{Step 5:} Since the automatic focusing of the lasers and the alignment of the plate lead to slight deviations from an ideal uniform grid, the $x$- and $z$-velocity fields are linearly interpolated to a smaller uniform grid. As Figure~\ref{fig:LDVoverview} shows, the uniform grid has almost the same grid width and is only used to correct the angles and reduce the number of measurement points. In addition, the uniform grid makes it possible to average the data in the $y$-direction, as mentioned above. This results in a preliminary array $\tilde \mV_{\Meas}$  of size $2 \times 11 \times 9 \times 782$, here the first dimension distinguishes between the $x$-and $z$-velocity, the second dimension is for the change in the $x$-coordinate, the third dimension is for the change in the $y$-coordinate, and the fourth dimension is for the change in time. 

\begin{figure}[ht!]
    \centering
    \subcaptionbox{$\bar v_x$ (in-plane) \\ $\SI{0.4}{\milli\meter}$ notch depth}
    {\includegraphics[height=.5\MyWidth, angle=-90]{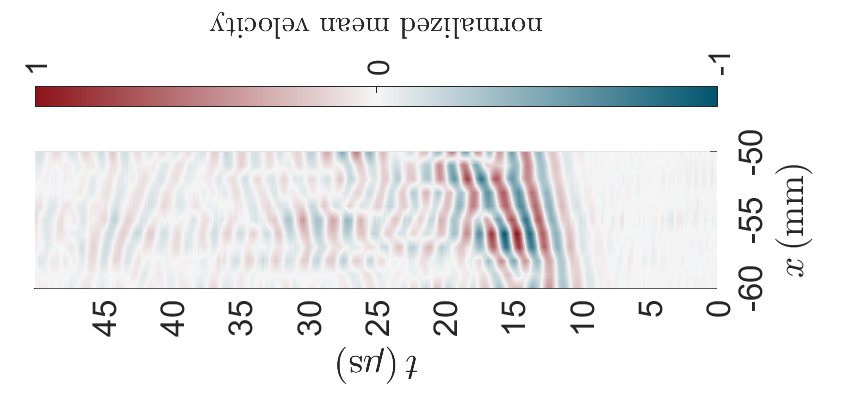}}
    \hfil
    \subcaptionbox{$\bar v_z$ (out-of-plane) \\ $\SI{0.4}{\milli\meter}$ notch depth}
    {\includegraphics[height=.5\MyWidth, angle=-90]{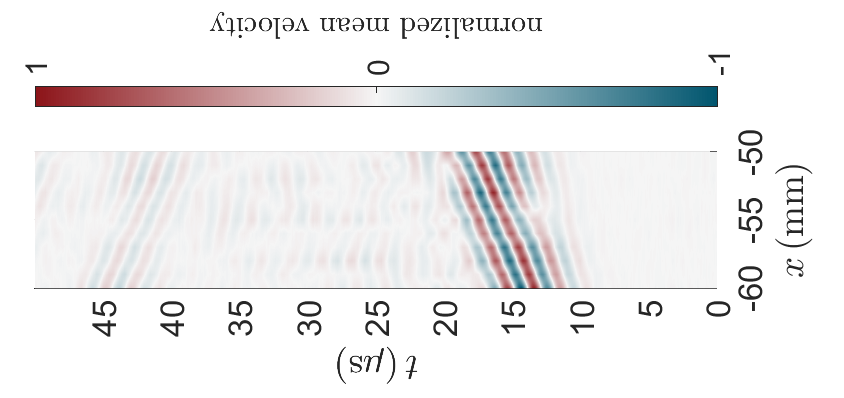}}
    \hfil
    \subcaptionbox{$\bar v_x$ (in-plane) \\ $\SI{0.8}{\milli\meter}$ notch depth}
    {\includegraphics[height=.5\MyWidth, angle=-90]{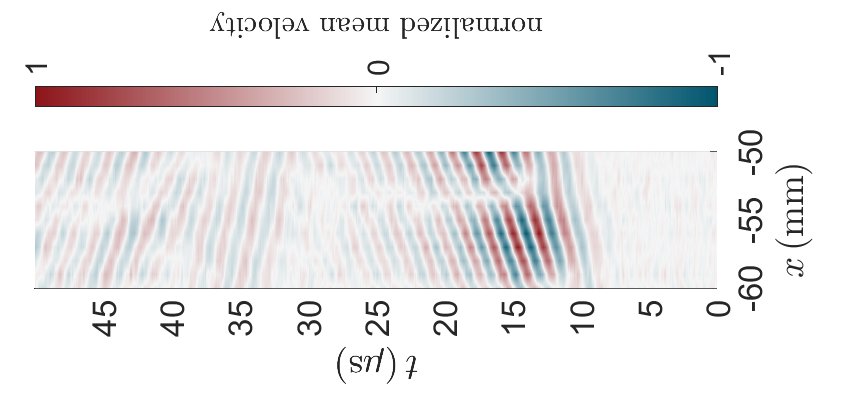}}
    \hfil
    \subcaptionbox{$\bar v_z$ (out-of-plane) \\ $\SI{0.8}{\milli\meter}$ notch depth}
    {\includegraphics[height=.5\MyWidth, angle=-90]{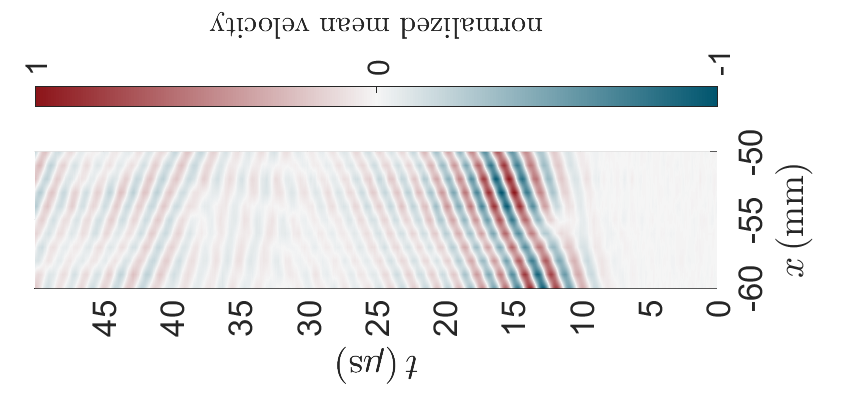}}
    \caption{Normalized mean velocity fields for the interpolated uniform grid (Figure~\ref{fig:LDVoverview}) stored in $\mV_{\Meas}$.}
    \label{fig:XTwavefield}
\end{figure}

Finally, the data is normalized to the largest amplitude.
In summary, this pre-processing of the preliminary array $\tilde \mV_{\Meas}$ is 
\begin{equation*}
    \mV_{\Meas} = \nrm \big( \mean_y (\tilde \mV_{\Meas} ) \big) \in \RR^{2 \times 11 \times 1 \times 782},
\end{equation*}
where
$\nrm(\mA) = \mA / \max\Big( \vek\big( \env_t(\mA) \big) \Big)$ 
is the normalization of an arbitrary array $\mA$ to its larges amplitude. 
The normalization is defined by an envelope operator over the time axis $\env_t(\mA)$ through a discrete Hilbert transform and a reshaping of a multidimensional array $\mA$ into a vector $\vek(\mA)$. To give an impression of the experimental data, Figure~\ref{fig:XTwavefield} shows the mean velocity fields stored in the array $\mV_{\Meas}$ for both plates. The data for both plates were recorded using the same steps. The only difference is that the original measuring grid is slightly different because the scanning area is defined by hand in the LVD software. In Figure~\ref{fig:XTwavefield}, mainly the time interval between $\SI{33}{\micro\second}$ and $\SI{50}{\micro\second}$ shows the reflected waves coming from the notches.

\section{The Inverse Problem}~\label{sec:Theory}


The defect reconstruction is formulated as an inverse problem. 
The geometric defect parameters $\vpara$ must be estimated so that the simulation model best fits the experimental data. The final reconstruction parameters are defined as
\begin{equation*}
    \vpara^{\min} = \arg \min_{\vpara \in \paraspace} \|\Fo(\vpara) - \vy_{\Meas}\|^2
\end{equation*}
with the forward operator $\Fo\colon \vpara \mapsto \vy_{\Sim} \in \RR^{\NMeas \times 1}$ (Subsection~\ref{sec:ForwardOperator}),
a vector $\vy_{\Meas} \in \RR^{\NMeas \times 1}$ containing the experimental data
and the parameter space $\paraspace$.

Our previous research has shown that similar inverse problems have significantly fewer local minima if the difference between the experimental data and the forward operator is based on the envelope of the signals~\cite{bulling2022defect}.
However, our previous studies evaluated a single displacement signal, while the current work considers the average velocity at multiple points.
Nevertheless, we came to the same conclusion after a similar study not presented here.
Consequently, the vector $\vy$ is defined by calculating the envelope over the time axis, and the array is reshaped into a vector, i.e.
\begin{equation}
    \vy_{\Meas} = \vek\big(\env_t(\mV_{\Meas})\big), \label{eq:Ymeas}
\end{equation}
where $\mV_{\Meas}$ is the multi-dimensional array described in Section~\ref{sec:Data}, and thus $\NMeas = 2\cdot11\cdot782 = \num{17204}$.

\subsection{The Forward Operator}\label{sec:ForwardOperator}
The forward operator mainly calculates the wave propagation in a cross-sectional model of a steel plate.
The main idea is to use a semi-analytical model based on the SBFEM, which leads to a comparatively small number of degrees of freedom in the final system of equations. The problem is solved in the (complex) frequency domain to obtain a model with very few degrees of freedom. The Exponential Window Method (EWM)~\cite{bulling2022defect,kausel2017advanced} is used to compute a finite domain in the frequency domain. 
The EWM is based on the Discrete Laplace Transformation $\DLT(\vs) = \FFT \big(\vs \odot \exp(-\EWM \vt) \big)$ and 
and its inverse $\DLT^{-1}(\hat \vs) = \FFT^{-1}(\hat \vs) \odot \exp(\EWM \vt)$, where 
$\FFT$ is the Fast Fourier Transform over the time dimension, 
$\vs$ and $\hat \vs$ are arbitrary signal arrays in the time and frequency domains, 
$\vt$ is the uniform time vector of the measurement, 
$\EWM=0.5\dOmega$ is the name-giving exponential factor,
$\dOmega$ is the resulting angular frequency step due to the sampling frequency of the time vector $\vt$,
$\exp$ is the entry-wise exponential function,
and $\odot$ is the entry-wise product along the time dimension.
In addition to using these transforms, the wave equation in the EWM is derived for discrete complex angular frequencies $\omega_\ell = (\ell \dOmega - \ima\EWM)\in\CC$.
Note that only the real part of the complex angular frequency changes with the index,
while the imaginary part remains constant.

\begin{figure}[ht!]
    \centering
    \begin{tikzpicture}[x=\MyUnit,y=\MyUnit]
    \node[inner sep=0pt] (Overview) at (0,0)
        {\includegraphics[width=9\MyUnit]{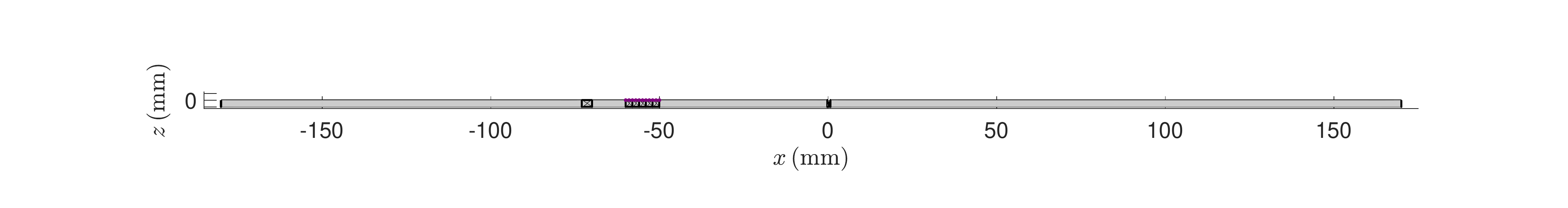}};
    \node[inner sep=0pt] (Detail11) at (-2.5,+1.5)
        {\fbox{\includegraphics[width=4\MyUnit]{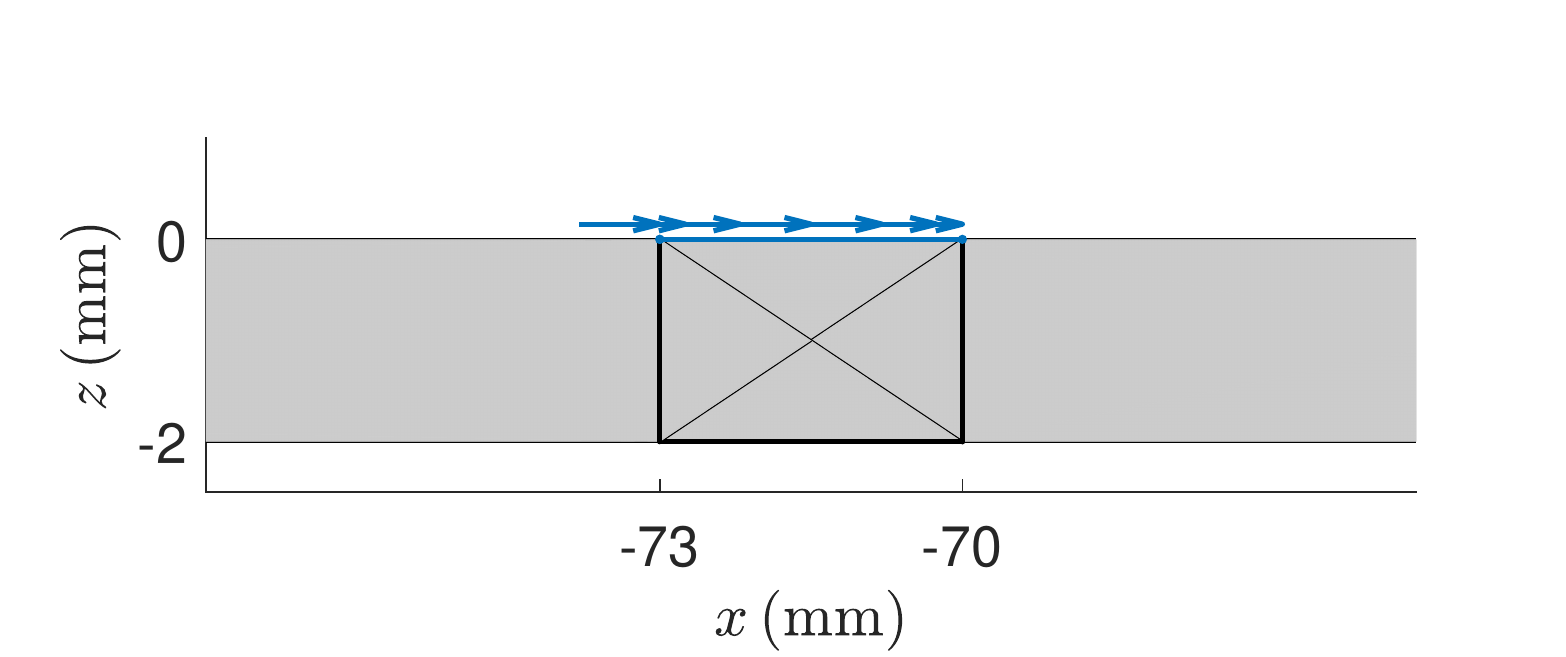}}};
    \draw[-,thick] (-1.15,0) -- (Detail11.south);
    \node[inner sep=0pt] (Detail12) at (+2.5,+1.5)
        {\fbox{\includegraphics[width=4\MyUnit]{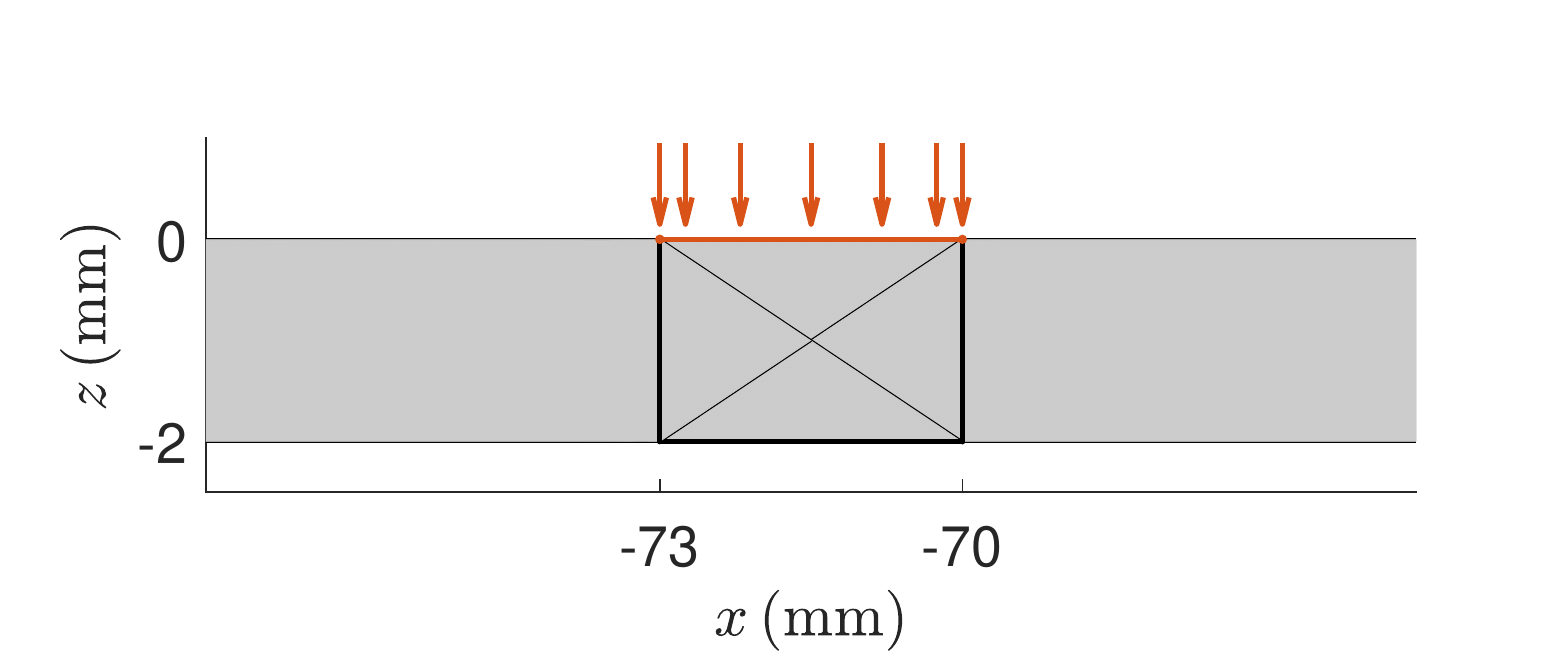}}};
    \draw[-,thick] (-1.15,0) -- (Detail12.south);
    \node[inner sep=0pt] (Detail21) at (-2.5,-1.5)
        {\fbox{\includegraphics[width=4\MyUnit]{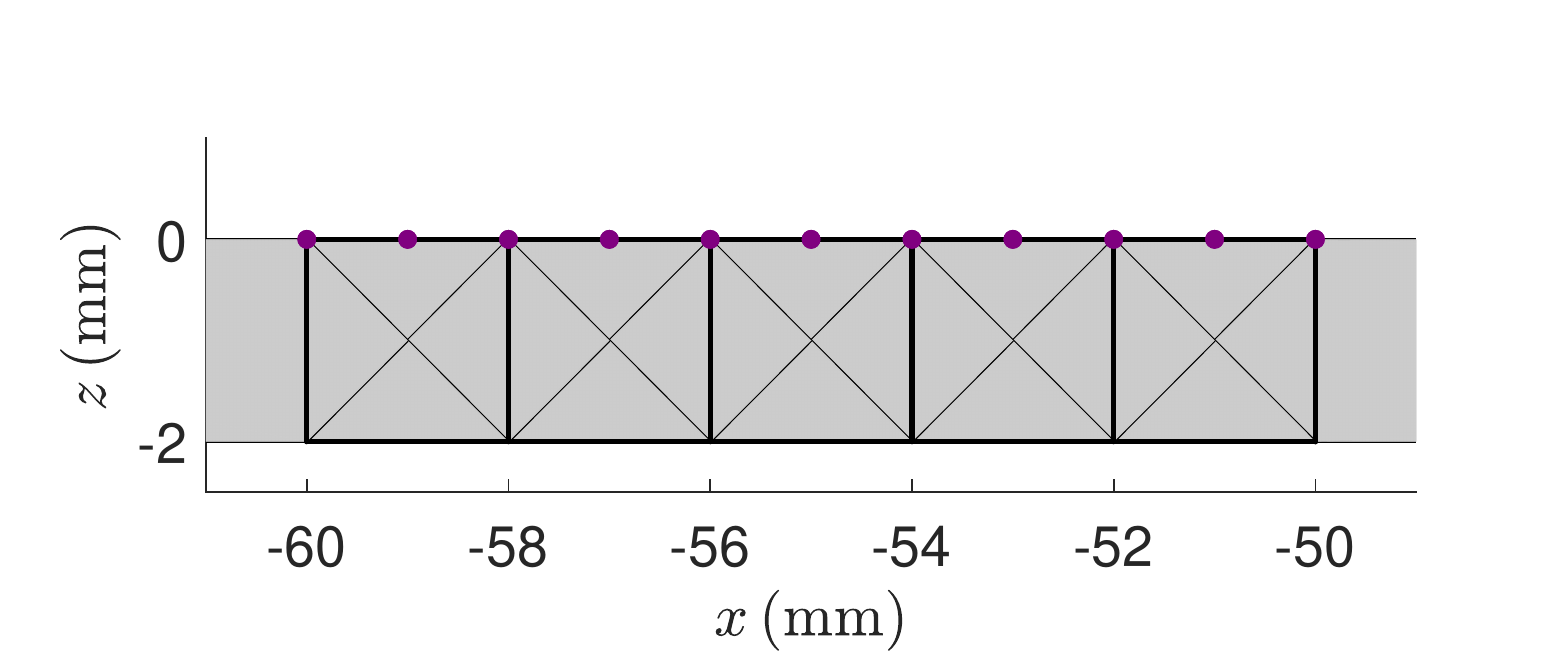}}};
    \draw[-,thick] (-.75,0) -- (Detail21.north);
    \node[inner sep=0pt] (Detail22) at (+2.5,-1.5)
        {\fbox{\includegraphics[width=4\MyUnit]{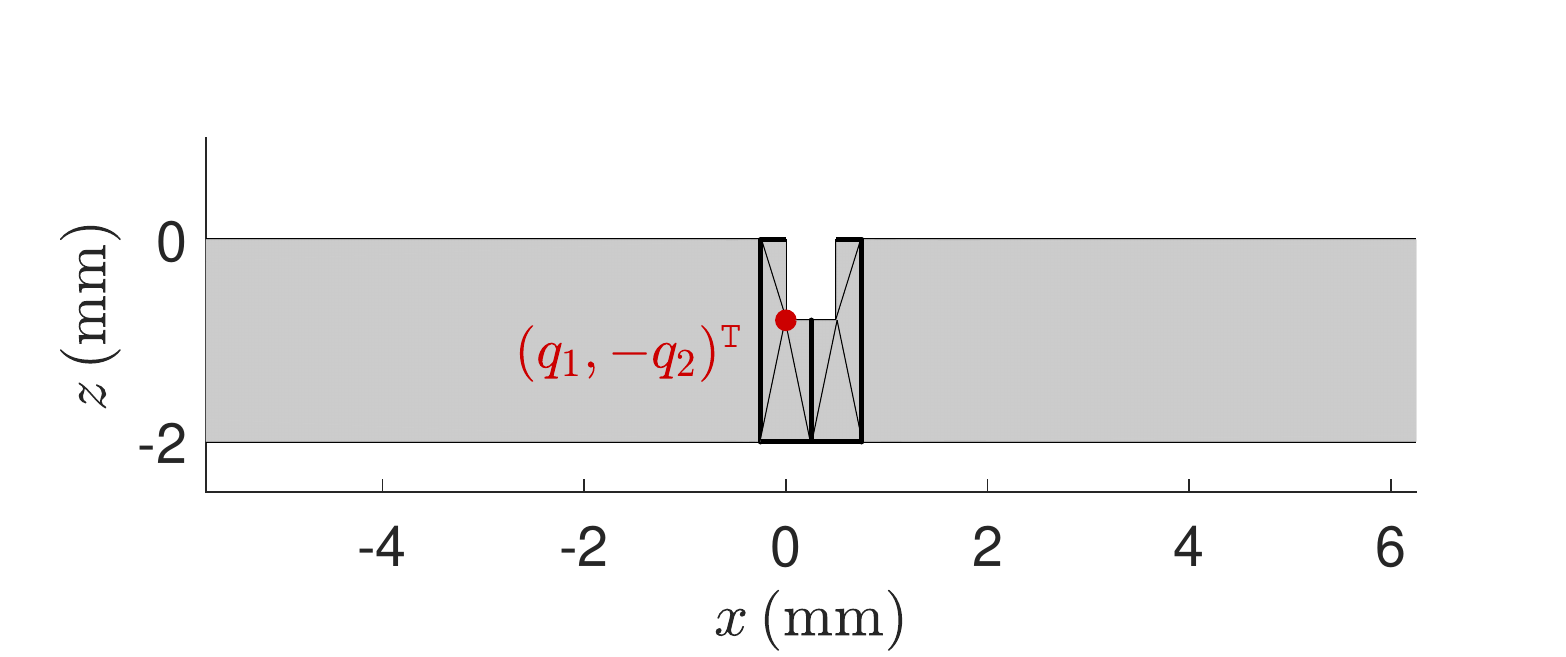}}};
    \draw[-,thick] (+0.25,0) -- (Detail22.north);
    \node[xshift=.1\MyUnit,yshift=0\MyUnit] at (-2.5,+2.15) {$\vbtau_1$};
    \node[xshift=.1\MyUnit,yshift=0\MyUnit] at (+2.5,+2.15) {$\vbtau_2$};
    \end{tikzpicture}
    \caption{Cross-sectional SBFEM model of the center line (Fig.~\ref{fig:LDVoverview}) of the plate with a $\vpara = (\SI{0}{\milli\meter},\SI{.8}{\milli\meter})^\tr$.\label{fig:SBFEMmodel}}
\end{figure}

In our and many other experimental setups, there is no validated piezoelectric model of the sensor used.
The current approach fits not only the damage to the experimental data but also the excitation model to account for this difficulty.
For this fitting, the experimental data is transformed into the frequency domain
\begin{equation}
    \hat \mV_{\Meas} = \DLT(\mV_{\Meas}),\label{eq:FreqData}
\end{equation}
where $\mV_{\Meas}$ is the multi-dimensional array described in Section~\ref{sec:Data}. 
The excitation model is based on two tractions acting on the plate boundary at the sensor region with unknown amplitudes and phases.
Therefore, the sensor itself is not modeled but only its force.
For each complex angular frequency $\omega_\ell$ in a certain frequency range with relevant amplitudes of the experimental data, the cross-sectional plate model (Figure~\ref{fig:SBFEMmodel}) approximates the two mean velocities $\fbv_j$, $j\in\{1,2\}$, one for each basic type of mean sensor traction
\begin{align}
    \vbtau_1 & = (1,0)^\tr && \text{and} & \vbtau_2 = (0,-1)^\tr \label{eq:Tractions}
\end{align}
and then, the two results $\fbv_j$ are weighted by a transfer function to fit the experimental data in the frequency domain (see Equations~\eqref{eq:Weighting}. The range between $\SI{10}{\kilo\hertz}$ and $\SI{1.5}{\mega\hertz}$ is selected as the relevant frequency range for the current application.

For simplicity, we omit the index dependence of $\ell$ for the complex angular frequency $\omega_\ell = \omega$ in the following.
For the cross-section at the center line (Figure~\ref{fig:LDVoverview}) the mean velocity fields are
$\fbv_j = (\bar{v}_j^x,\bar{v}_j^z)^\tr = (\ima \omega)\fbu_j$ are based on the elastic wave equations
\begin{align}
    (\ima \omega)^2 \rho \fbu_j(\vx) & = \vnabla \cdot \vsigma\left(\fbu_j(\vx)\right) & \vx & \in \Omega, \label{eq:PDEa}\\
    \vsigma\left(\fbu_j(\vx)\right) \fn & = \vbtau_j(\vx) & \vx & \in \Gamma_{\text{sensor}}, \label{eq:PDEb}\\
    \vsigma\left(\fbu_j(\vx)\right) \fn & = \vo & \vx & \in \Gamma_{\text{free}}, \label{eq:PDEc}
\end{align} 
where 
$\rho$ the density,
$\fbu_j = (\bar{u}_j^x,\bar{u}_j^z)^\tr$ the mean displacement,
$\vsigma$ the mean linear stress operator,    
$\fn$ the outer normal vector, 
$\vbtau_j$ the mean traction,
and the cross-section coordinates $\vx = (x,z)^\tr$.
The linear stress tensor $\vsigma$ is defined by Hooke's law, i.e, $\vsigma(\fbu_j) = \mC:\vepsi(\fbu_j)$, and the linear strain-displacement relation, $\vepsi(\fbu_j) = 0.5\left((\vnabla\fbu_j)^\tr + (\vnabla\fbu_j)\right)$, where $\mC$ is the fourth-order elasticity tensor with the plane strain assumption given by the material parameters.
The material parameters for the plates studied are density $\rho = \SI{7.9}{\gram\per\cubic\centi\meter}$, Young's modulus $\EE = \SI{200}{\giga\pascal}$, and Poisson's ratio $\nu = \num{0.3}$.

The cross-section at the center line is the domain $\Omega$ with a parameterized notch, which can be expressed mathematically as
\begin{equation*}
    \Omega = (\SI{-180}{\milli\meter},\SI{170}{\milli\meter})\times(\SI{0}{\milli\meter},\SI{-2}{\milli\meter})
    \setminus
    \overline{(q_1,q_1+0.5\si{\milli\meter})\times(\SI{0}{\milli\meter},- q_2)},
\end{equation*}
where $q_1$ is the $x$-position of the notch,
and $q_2$ is the depth of the notch.
The tractions are applied to the sensor area $\Gamma_{\text{sensor}} = [\SI{-73}{\milli\meter},\SI{-70}{\milli\meter}]\times \{0\}$, while the rest of the boundary $\Gamma_{\text{free}} = \partial\Omega \setminus \Gamma_{\text{sensor}}$ is traction-free.

To apply the SBFEM to the forward model  (Figure~\ref{fig:SBFEMmodel}), the domain is subdivided into super-elements. For the current forward model, two different types of super-elements are used, which are briefly summarized below.
\begin{figure}[ht!]
    \centering
    \subcaptionbox{Waveguide super-element.\label{fig:SBFEMsEleA}}
    {\includegraphics[width=\MyWidth]{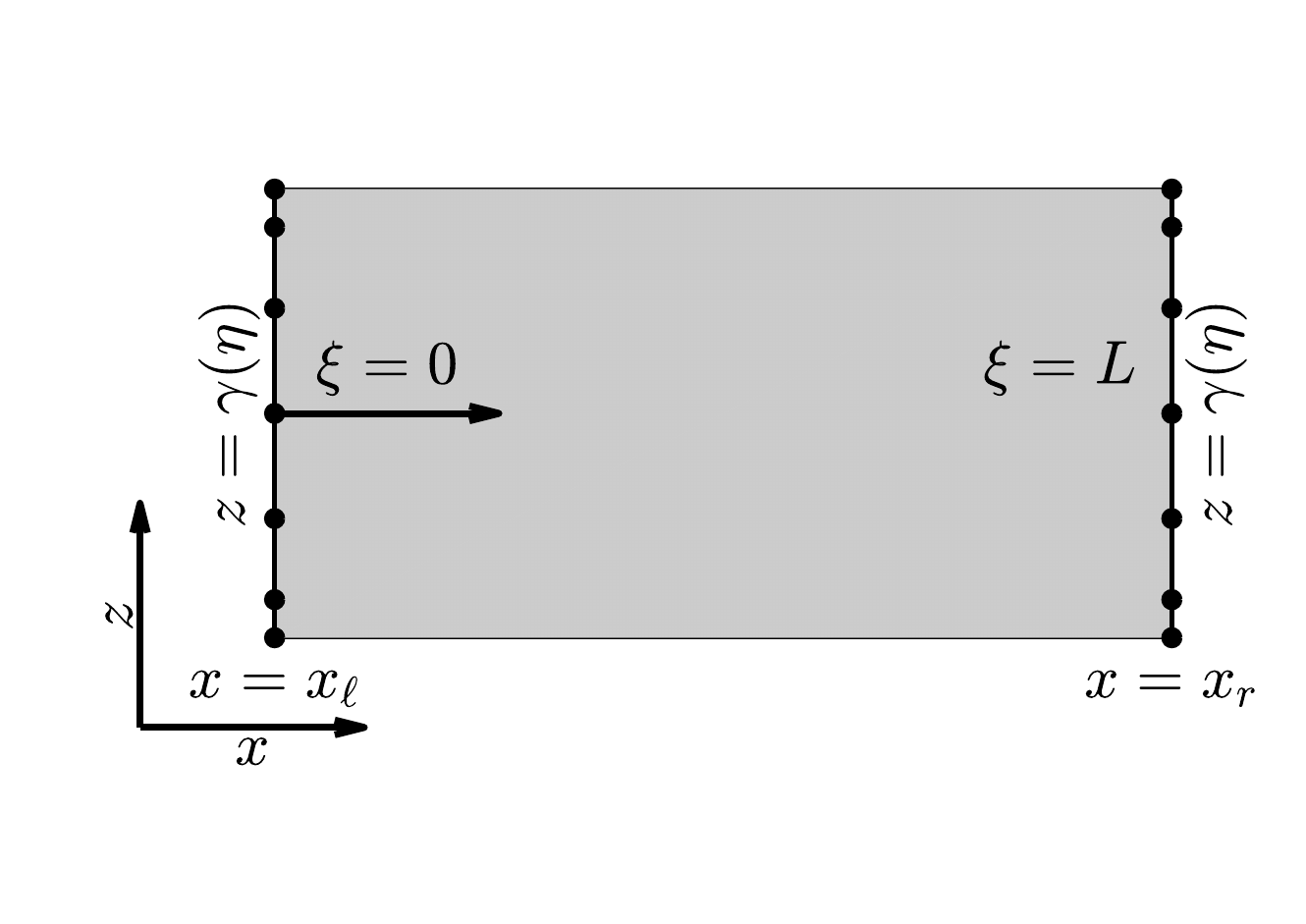}}
    \hfil
    \subcaptionbox{Polygonal super-element.\label{fig:SBFEMsEleB}}
    {\includegraphics[width=\MyWidth]{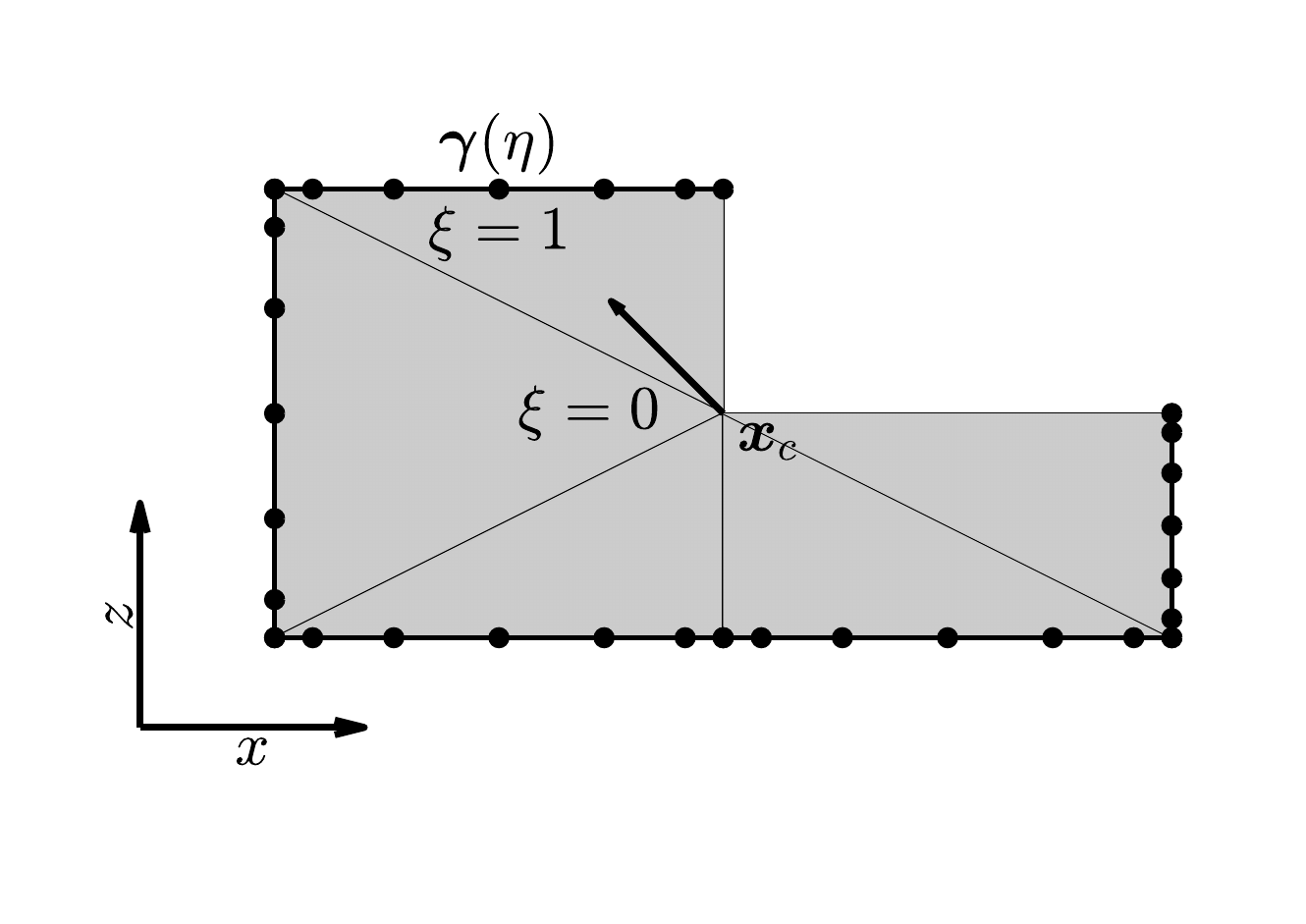}}
    \caption{The two types of SBFEM super-elements.}
\end{figure}

The first type of super-element approximates the displacement in the undisturbed parts of the plate's cross-section. Gravenkamp et al. developed this super-element for waveguides~\cite{gravenkamp2015simulation}; therefore, it will be referred to as waveguide super-element in the following.
Figure~\ref{fig:SBFEMsEleA} shows an example of a super-element, where $\gamma(\eta)$ is a 1D finite element mesh of the left and right edges in the $z$-direction with an element-wise parameterization coordinate $\eta$.
The coordinates $x_\ell$ and $x_r$ are the position of the left and right boundaries. Using the coordinate transformation
\begin{equation*}
    \begin{pmatrix}
        x \\
        z
    \end{pmatrix}
    =
    \begin{pmatrix}
        \xi + x_\ell \\
        \gamma(\eta)
    \end{pmatrix},
\end{equation*}
with $\xi$ between $0$ and the super-element length $L = x_r - x_\ell$,
a product ansatz
\begin{align}
    \fbu(x,z) = \fN(\eta) \fub(\xi),\label{eq:Ansatz}
\end{align}
is assumed, where $\fN(\eta)$ is the matrix of the vector-valued shape functions defined by the finite element mesh and $\fub(\xi)$ is the vector of nodal displacements~\cite{gravenkamp2015simulation,gravenkamp2018efficient}.
While in general a variety of different shape functions can be adopted~\cite{gravenkamp2021high}, for this work, spectral shape functions of degree $\pp=6$ (with 7 nodes) based on the Gauss-Lobatto integration points are used. Substituting the product ansatz~\eqref{eq:Ansatz} into the Equations~\eqref{eq:PDEa} and integrating by parts, we obtain the Ordinary Differential Equation (ODE) called SBFEM equation in displacement 
\begin{align}
    \mE_0 \partial_{\xi\xi} \fub(\xi) + (\mE_1^\tr - \mE_1) \partial_\xi \fub(\xi) - \mE_2 \fub(\xi) + \omega^2 \mM_0 \fub(\xi) & = \vo, \label{eq:SBFEM_u_WG}
\end{align}
where the coefficient matrices $\mE_i$ and $\mM_0$ are calculated similarly to the finite element method by integrating the shape functions and their derivatives. 
The ODE~\eqref{eq:SBFEM_u_WG} leads to an quadratic eigenvalue problem
\begin{align*}
    \mE_0 (\ima \wn)^2 \vpsi + (\mE_1^\tr - \mE_1) (\ima \wn) \vpsi - \mE_2 \vpsi + \omega^2 \mM_0 \vpsi & = \vo,
\end{align*}
where the eigenvalue $(\ima \wn)$ is the product between the wavenumber of a mode and the imaginary unit, and the eigenvector $\vpsi$ is the nodal mode shape vector.
The solution for the nodal displacement has the form
\begin{equation*}
    \fub(\xi) = 
      \mPsi_\ell\diag\Big(\exp\big(\ima \vk_\ell \xi \big)\Big)\hat \vc_\ell 
    + \mPsi_r   \diag\Big(\exp\big(\ima \vk_r(\xi-L) \big)\Big)\hat \vc_r,
\end{equation*}
where 
$\mPsi_\ell$, $\mPsi_r$ are the eigenvector matrices with columns $\vpsi$,
$\hat \vk_\ell$, $\hat \vk_r$ are vector with entries $\wn$,
$\mA = \diag(\va)$ is the diagonal matrix with entries according to the vector $\va$.
The coefficient vectors $\hat \vc_\ell$, $\hat \vc_r$ are the unknowns of the wave field for this super element and must be determined by the boundary values. 
The index $\ell$ indicates the waves propagating from left to right, while the index $r$ indicates the waves propagating from right to left.
A rearrangement of the equation for the displacements $\fub(\xi)$ and an equation for the internal nodal forces $\left(\mE_0\, \partial_{\xi} + \mE_1^\tr\right) \fub(\xi)$ with $\xi = 0$ and $\xi = L$ then allows to determine the local stiffness matrix $\hat \mS^{b}$ for this super-element. For a detailed derivation, the interested reader is referred to \cite{gravenkamp2015simulation,gravenkamp2018efficient}. We want to emphasize that the length $L$ of the waveguide super-element does not affect its computational cost, making it an ideal choice for the long, undamaged parts of the cross-section. 

The second type of super-element approximates the displacement for a star-convex polygonal sub-domain. Figure~\ref{fig:SBFEMsEleB} shows a 1D finite element mesh $\vgamma$ and a special point $\SC$ called the scaling center. Similar to the first type of super-element, a coordinate transformation
\begin{equation*}
    \vx = \xi\big(\vgamma(\eta) - \SC\big) + \SC
\end{equation*}
with $\xi$ between $0$ and $1$,
the product ansatz form Equation~\eqref{eq:Ansatz} and integration by parts are utilized to write the partial differential Equation~\eqref{eq:PDEa} into another SBFEM-equation in displacement
\begin{align}
    \xi^2 \mE_0 \partial_{\xi\xi} \fub(\xi) + \xi (\mE_0 + \mE_1^\tr - \mE_1) \partial_\xi \fub(\xi) - \mE_2 \fub(\xi) + \xi^2 \omega^2 \mM_0 \fub(\xi) & = \vo. \label{eq:SBFEM_u_Poly}
\end{align}
Since there is no known closed-form solution for this type of ODE~\eqref{eq:SBFEM_u_Poly}, the approximation approach for polygonal super-elements is more sophisticated. A differential equation for the local stiffness matrix $\hat \mS^{b}$ is derived and approximated by a continued-fraction iteration.
Defining the local stiffness as a matrix expression for the inner forces, i.e.,
\begin{align*}
	\hat \mS^{b}(\xi,\omega)\, \fub(\xi) = \left(\xi\,\mE_0\, \partial_{\xi} + \mE_1^\tr\right) \fub(\xi),
\end{align*}
considerations regarding the scalability of the wave equation for linear elastic problems lead to 
the assumption that the stiffness matrix is a function of the single quantity $\xo = (\ima \xi \omega)^2$, i.e., $\hat \mS^{b}(\xi,\omega) = \hat \mS^{b}(\xo)$. 
Substituting the these assumptions into Equation~\eqref{eq:SBFEM_u_Poly} leads to
\begin{align}
	2 \xo\, \partial_\xo \hat \mS^b(\xo) + (\hat \mS^b(\xo) - \mE_1)\mE_0^{-1}(\hat \mS^b(\xo) - \mE_1^{\tr}) - \mE_2 - \xo \mM_0 & = \vo.\label{eq:SBFEM_S_poly}
\end{align}
Subsequently, Equation~\eqref{eq:SBFEM_S_poly} is approximated by a continued-fraction iteration of order $M$
\begin{align*}
    \hat \mS^{b}(\xo)  
    & = \hat \mS^{(0)}(\xo), \\
    \hat \mS^{(m)}(\xo)
    & = \mS_0^{(m)} + \xo \mS_1^{(m)}
    - \xo^2 \mX^{(m+1)} \big(\hat \mS^{(m+1)}(\xo)\big)^{-1}  (\mX^{(m+1)})^\tr, \\
    \hat \mS^{(M)}(\xo) & = \mS_0^{(M)} + \xo \mS_1^{(M)},
\end{align*}
where
$\mS_0^{(m)}$, $\mS_1^{(m)}$ are the continued-fraction matrices  of the $m$-th term,
and $\mX^{(m)}$ is the $m$-th pre-conditioner matrix. The details of the computation of the matrices $\mS_0^{(m)}$, $\mS_1^{(m)}$ and $\mX^{(m)}$ can be found in \cite{chen2014high}. The matrices are derived either by algebraic Riccati or Lyapunov equations. Both types of algebraic equations can be solved by an eigenvalue problem~\cite{bulling2022defect,song1997scaled}. Generally, the order $M$ of the iteration can be adjusted, but we followed the simple rule $M=\pp$, where $\pp=6$ is the polynomial degree of the shape functions. 
It is worth noting that the displacement in the re-entrant corner, which leads to reduced convergence in many numerical methods, is well approximated by the discretization used here. If a re-entrant corner is located in the scaling center, optimal convergence is still observed in SBFEM~\cite{bulling2019comparison}. This property is the reason for the two polygonal super elements around the notch (see Figure~\ref{fig:SBFEMmodel}).

\begin{algorithm}
\caption{Forward Operator}
\begin{algorithmic}[1]
  \Statex \textbf{Input:} geometric parameter $\vpara$, time vector $\vt$, the measurement array $\mV_{\Meas}$, regularization parameter $\beta$ 
  \vspace{.5ex} \hrule \vspace{.5ex}
  \Statex {\color{gray}\textit{Block of all frequency-independent quantities:}}
  \State compute $\zeta$ from $\vt$
  \State $\hat \mV_{\Meas} = \DLT(\mV_{\Meas})$
  \State allocate $\hat \mV_{\Sim}$, $\mA$, $\mF$
  \State compute $\mA$
  \For{all finite elements $e$}
    \State compute $\mF^{e}$
    \State add $\mF^{e}$ to $\mF$
  \EndFor
  \For{all super elements}
    \State compute and store $\mM_0$,$\mE_0$,$\mE_1$,$\mE_2$,$\mS_0^{(m)}$,$\mS_1^{(m)}$,$\mX^{(m+1)}$ depending on the type of the super element
  \EndFor
  \Statex 
  \Statex {\color{gray}\textit{Block of all frequency-dependent quantities:}}
  \State determine a relevant frequency range $\vomega$ from $\hat \mV_{\Meas}$
  \For{all $\omega$ in $\vomega$}
    \State allocate $\hat \mS$
    \For{all super elements}
      \State compute $\hat \mS^{b}$ depending on the type of the super element
      \State add $\hat \mS^{b}$ to $\hat \mS$
    \EndFor
    \State $\hat \mU = \hat\mS \backslash \mF$
    \State $\hat \mV^{}_{\Sim} = (\ima \omega)\mA\hat \mU$
    \State $\hat \vv_{\Meas} = \vek\big( \hat \mV_{\Meas} \big|_{\omega} \big)$
    \State $\hat \vh =(\hat \mV^\ct_{\Sim} \hat 
    \mV^{}_{\Sim} + \beta \mI) \backslash (\hat \mV_{\Sim}^\ct\hat \vv^{}_{\Meas})$
    \State $\hat \vv^{}_{\Sim} = \hat \mV^{}_{\Sim} \hat \vh$
    \State sort $\hat \vv^{}_{\Sim}$ into $\hat \mV_{\Sim}$
  \EndFor
  \Statex
  \Statex {\color{gray}\textit{Transform back to time domain and final steps:}}
  \State $\mV_{\Sim} = \DLT^{-1}(\hat \mV_{\Sim})$
  \State $\vy_{\Sim} = \vek\big(\env_t(\mV_{\Sim})\big)$
\end{algorithmic}
\label{Algo:Forward}
\end{algorithm}

The local stiffness matrices $\hat \mS^{b}$ are assembled to a global stiffness matrix $\hat \mS$. To approximate the equations~\eqref{eq:PDEa}-\eqref{eq:PDEc}, the two traction terms ( Equation~\eqref{eq:Tractions}) must also be taken into account. Since the boundary on which the forces act is discretized with finite boundary elements, the nodal versions of the tractions $\mF$ are represented analogously to the FEM by a $\NDoF \times 2$ matrix, where $\NDoF$ is the number of degrees of freedom, and each column is associated with a traction term in Equation~\eqref{eq:Tractions}. Note that the nodal tractions $\mF$ are frequency independent.
A linear system of equations must be solved for each frequency
\begin{align*}
    \hat\mS \hat \mU & = \mF
    && \text{ with } \hat \mS \in \CC^{\NDoF \times \NDoF} \text{ and }\hat \mU,\mF \in \CC^{\NDoF \times 2}.
\end{align*}
For the model in Figure~\ref{fig:SBFEMmodel}, the number of degrees of freedom  $\NDoF$ is $\num{346}$.

As mentioned above, the concept of a transfer function is utilized to generate the simulated velocity.
Since the simulation model is constructed such that there are nodes at the measurement points, the velocities $\hat \mV^{}_{\Sim}$ for the two traction terms (Equation~\eqref{eq:Tractions}) are calculated by
$\hat \mV^{}_{\Sim} = (\ima \omega)\mA\hat \mU$ with the assignment matrix
\begin{equation*}
    (\mA)_{mn} =
    \begin{cases}
        1 & \text{if the $n$-th simulation degree corresponds to the $m$-th measurement degree}\\
         0          & \text{otherwise}
    \end{cases}
    \in \RR^{\NEval\times\NDoF}.
\end{equation*}
The two columns of velocities are afterward fitted to the measurement data
\begin{equation*}
    \hat \vv_{\Meas} = \vek\big( \hat \mV_{\Meas} \big|_{\omega} \big) \in \CC^{\NEval \times 1},
\end{equation*}
where $\hat \mV_{\Meas} \big|_{\omega}$ is the part of the array in Equation~\eqref{eq:FreqData} for the current frequency $\omega$ with $\NEval = 2 \cdot 11 = 22$, by the transfer functions
\begin{equation}
    \hat \vh = (\hat \mV^\ct_{\Sim} \hat \mV^{}_{\Sim} + \beta \mI)^{-1} \hat \mV_{\Sim}^\ct\hat \vv^{}_{\Meas}.\label{eq:TransferFunctions}
\end{equation}
A small regularization parameter $\beta = \num[output-exponent-marker = \ensuremath{\mathrm{E}}]{1E-9}$ is introduced, because the invertibility of $\hat \mV^\ct_{\Sim} \hat \mV^{}_{\Sim}$ is not guaranteed.
The final simulated velocities for the current frequency at the measurement points is
\begin{equation}
    \hat \vv^{}_{\Sim} = \hat \mV^{}_{\Sim} \hat \vh. \label{eq:Weighting}
\end{equation}
For all relevant frequencies, the simulation results are collected, reshaped and transferred back into the time domain by the inverse Laplace transformation $\DLT^{-1}$ in an array $\mV_{\Sim}$.
Analogously to Equation~\eqref{eq:Ymeas}, the final output of the forward operator is
\begin{equation*}
    \vy_{\Sim} = \vek\big(\env_t(\mV_{\Sim})\big).
\end{equation*}

A simplified version of the forward operator is shown in Algorithm~\ref{Algo:Forward}. This version is intended as a guide for a first implementation. However, many steps can be taken to speed up this version. For example, quantities that are independent of the geometric parameters $\vpara$, such as $\mF$, $\mA$, $\vomega$, $\hat \mV_{\Meas}$, can be calculated only once and treated as additional input values. In our implementation, the frequency-independent matrices, e.g., $\mM_0$,$\mE_0$,$\mE_1$,$\mE_2$, are stored for all super-elements that are independent of the geometric parameters $\vpara$. These super-elements are to the left of $x=\SI{-50}{\milli\meter}$.

\subsection{The Gradient Descent and the Initial Value}~\label{sec:IRGNM}
The reconstruction algorithm utilizes gradient-based optimization strategies. 
In contrast to typical derivative-free optimization methods gradient-based optimization methods have provable convergence properties such as guaranteed local convergence with high convergence speed.
However, such algorithms determine only local minima. Since the solution of the parameter estimation problem considered here is given by the global minimizer, for gradient-based methods it is important to provide a good initial guess which is in proximity of the globally optimal point.
We achieve this goal by computing the error between simulation and measurement at 100 random points uniformly distributed throughout the parameter space $\paraspace$.
The point with the smallest error is selected as an initial guess for the optimization procedure.

As an optimization approach, we selected a regularized Gauss-Newton method \cite{irgnm}, see Algorithm~\ref{Algo:IRGNM}.
\begin{algorithm}[t]
\caption{Iteratively Regularized Gauss-Newton Method (IRGNM)}
\begin{algorithmic}[1]
 \Statex \textbf{Input:} initial value $\vpara_0$, regularization parameter sequence $(\alpha_n)_{n\in {\NN}}\rightarrow 0$, maximal number of iterations $N$
 \vspace{.5ex} \hrule \vspace{.5ex}
  \For{$n = 0,\dots, N$}
 \State $\vpara_{n+1} = \vpara_n + (\Fo'(\vpara_n)^\ct\Fo'(\vpara_n)+ \alpha_n \mI)^{-1} (\Fo'(\vpara_n)^\ct (\vy_{meas}-\Fo(\vpara_n)) + \alpha_n (\vpara_0-\vpara_n) )$
 \If{$\|\vpara_{n+1}-\vpara_n\|<\epsilon$,}
    \State $\vpara^{\min} = \vpara_{n+1}$
    \State STOP
 \EndIf
 \EndFor
\end{algorithmic}
\label{Algo:IRGNM}
\end{algorithm}
In each iteration $n$ the iteratively regularized Gauss-Newton method minimizes a linearization of the original nonlinear least-squares problem (the error between simulation and measurement) with an additional Thikonov regularization term that vanishes as $n\rightarrow \infty$. 
The method is well understood and analyzed and is predominantly used in ill-posed parameter identification problems. 
Classical theory also guarantees convergence results in the presence of measurement errors.
Different than in prior cases \cite{bulling2022defect} when using synthetic data instead of real measurement data it was not necessary to utilize regularization. 
Hence, we set $\alpha_n=0$.

This method is highly dependent on accurate and cheaply available gradient information. Derivatives could be computed numerically by applying a finite difference scheme to the forward operator.
However, finite difference approximations of the derivative are usually very error prone and may be very costly to compute if the function to be differentiated has many variables. 
In the past we have experienced some issues with accuracy in this context and thus chose to use a different approach: 
As the whole simulation code is available to us we apply \textit{Algorithmic Differentiation} (AD) \cite{GrWa08,Nau12}. 
We selected the AD tool AdiMat \cite{adimat} as it can differentiate  simulations written in MATLAB. Furthermore, AdiMat can handle some technical programming issues such as datastructures etc. that other AD tools for MATLAB can not.
Overall an examination of the computed derivative values showed that the derivatives were within negligible tolerance of expected derivative values.

\section{Evaluation of the Proposed Algorithm}~\label{sec:Results}

In this section, we present the performance of our reconstruction algorithm for the two studied plates with different notch depths. The reconstruction is performed for a single parameter space that includes both notch depths. 
The position of the left notch $q_1$ is in the range between $\SI{-40}{\milli\meter}$ and $\SI{+40}{\milli\meter}$, which means that the distance is more than 17 A0-wavelengths.
The notch depth is between $\SI{0.1}{\milli\meter}$ and $\SI{1.1}{\milli\meter}$. The notch width is set to $\SI{0.5}{\milli\meter}$, the nominal width of both plates. The notch width is not of much interest because the notch is only used as a replacement geometry for a crack. Normally, cracks are assumed to be infinitesimally wide, so no width needs to be reconstructed.
It is essential that the numerical model is stable for the entire parameter space, since the initial values of the experiments are randomly drawn from this parameter space, as described in Section~\ref{sec:IRGNM}. 

\begin{figure}[ht!]
    \centering
    \subcaptionbox{Plate with a notch of $\SI{0.4}{\milli\meter}$ depth.}
    {\includegraphics[height=\MyWidth, angle=-90]{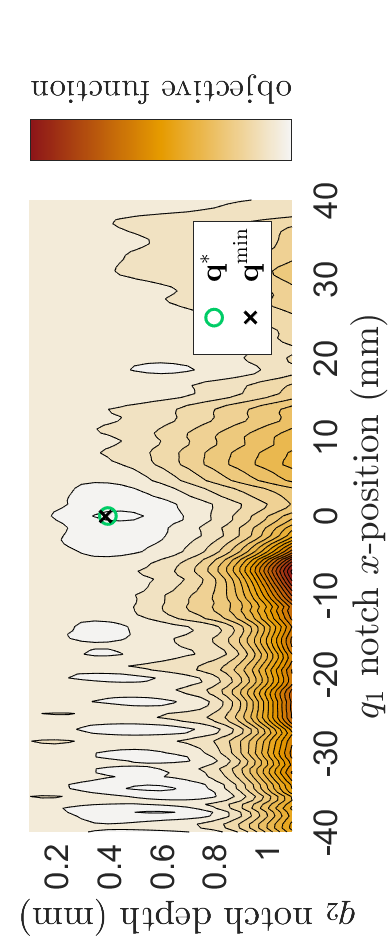}}
    \hfil
    \subcaptionbox{Plate with a notch of $\SI{0.8}{\milli\meter}$ depth.}
    {\includegraphics[height=\MyWidth, angle=-90]{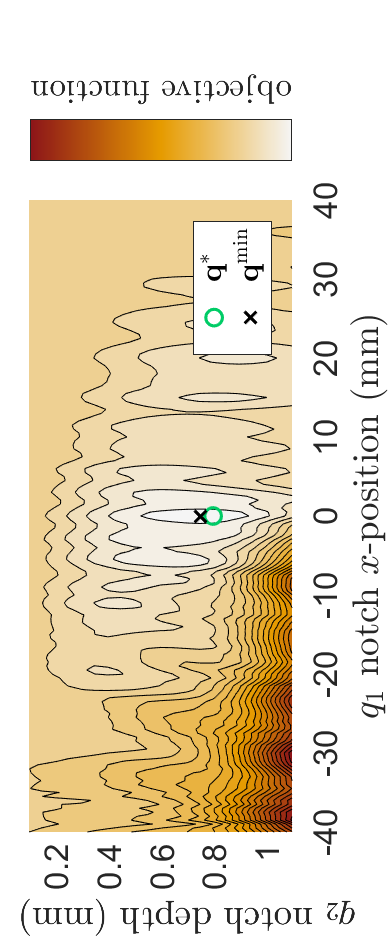}}
    \caption{Contour plot of the objective function with the nominal parameters $\vpara^*$ and the reconstruction parameters $\vpara^{\min}$.}
    \label{fig:ObjectiveFunction}
\end{figure}

As a preliminary investigation to the optimization, the objective function is examined. This additional investigation only provides insight into the optimization problem but is not necessary for the actual reconstruction. 
Figure~\ref{fig:ObjectiveFunction} shows the objective function for both plates over the entire examined parameter space.
The colored contour plots are based on a parameter grid of size $161 \times 21$. Both plates show a complex shape for the objective function with several local minima. However, the global minimum is located near the nominal parameters $\vpara^* = (\SI{0}{\milli\meter},\SI{.4}{\milli\meter})^\tr$ and $\vpara^* = (\SI{0}{\milli\meter},\SI{.8}{\milli\meter})^\tr$, respectively. There are several reasons why the global minimum of the objective function is not exactly at the nominal parameters.
First, the measurement naturally contains a certain amount of noise. Also, the assumptions for the numerical model may not be accurate enough. It should be emphasized that this is only a 2D cross-sectional model.
In addition, the plates were machined with a certain tolerance so that the actual value may differ from the nominal parameters.
In particular, the MFC sensors were manually aligned to the predefined position and glued on by hand. 
Furthermore, Figure~\ref{fig:ObjectiveFunction} shows that the reconstruction parameters in the area with a low objective function value are close to the nominal parameters.

\begin{table}[ht!]
    \begin{strech}
        \centering
        \begin{tabular}{c|c|c}
            description &
            nominal parameters&
            reconstruction parameters \\ 
            \hline
            plate with $\SI{.4}{\milli\meter}$ notch depth&
            $\vpara^{*} = (\SI{0}{\milli\meter},\SI{.4}{\milli\meter})^\tr$&
            $\vpara^{\min} = (\SI{-0.05}{\milli\meter},\SI{+0.39}{\milli\meter})^\tr$\\
            plate with $\SI{.8}{\milli\meter}$ notch depth&
            $\vpara^{*} = (\SI{0}{\milli\meter},\SI{.8}{\milli\meter})^\tr$&
            $\vpara^{\min} = (\SI{-0.11}{\milli\meter},\SI{+0.75}{\milli\meter})^\tr$\\
        \end{tabular}
        \caption{Reconstruction parameters}
        \label{tab:Parameters}
    \end{strech}
\end{table}

\begin{figure}[ht!]
    \centering
    \subcaptionbox{Plate with a notch of $\SI{0.4}{\milli\meter}$ depth.\label{fig:TransferFunctionsA}}
    {\includegraphics[width=\MyWidth]{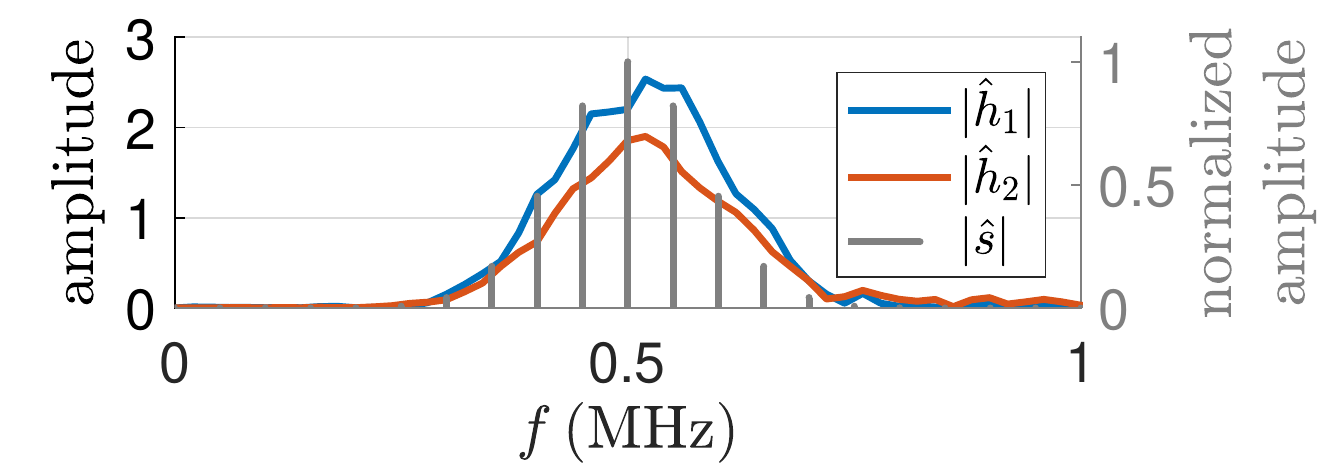}}
    \hfil
    \subcaptionbox{Plate with a notch of $\SI{0.8}{\milli\meter}$ depth.\label{fig:TransferFunctionsB}}
    {\includegraphics[width=\MyWidth]{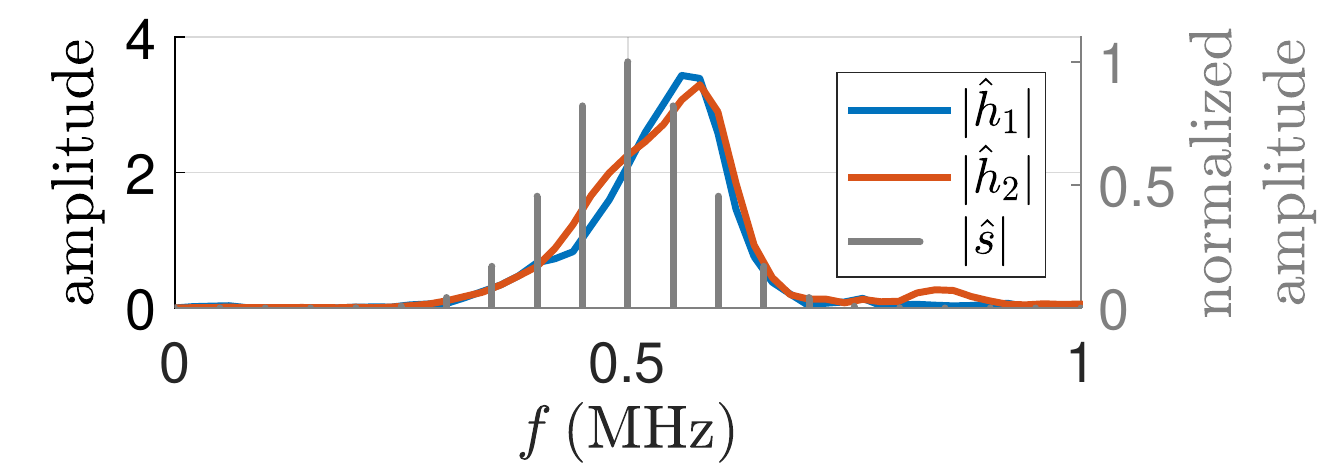}}
    \caption{Transfer functions of the model for the reconstruction parameter $\vpara^{\min}$ in the frequency domain. Colors correspond to those in Figure~\ref{fig:SBFEMmodel}.}
    \label{fig:TransferFunctions}
\end{figure}

The final reconstruction with the corresponding parameters $\vpara^{\min}$ (Table~\ref{tab:Parameters}) will be examined in the following. Figure~\ref{fig:TransferFunctions} shows the spectrum of the calculated transfer functions $\hat{h}_i$ (Equation~\eqref{eq:TransferFunctions}). It can be seen that both tractions make a significant contribution to the signal. Moreover, the maximum value of the transfer functions is close to the center frequency of $\SI{500}{\kilo\hertz}$ of the input signal $\hat{s}$. However, there are also qualitative differences between the two cases. While Figure~\ref{fig:TransferFunctionsA} shows an excellent agreement in the overall shape of the spectra between the transfer functions $\hat{h}_i$ and the input signal $\hat{s}$, a clear difference in shape can be observed in Figure~\ref{fig:TransferFunctionsB}.
The transfer functions represent not only the different MFC sensors, but also the adhesive bond to the plate. The fact that the adhesive bond is handmade could explain the differences, despite the same design of the sensors and exactly the same signal generator and amplifier. However, the difference shows that calculating the transfer functions makes the optimization more robust to variations in sensor performance.

\begin{figure}[ht!]
    \centering
    \subcaptionbox{Plate with a notch of $\SI{0.4}{\milli\meter}$ depth.}
    {\includegraphics[width=\MyWidth]{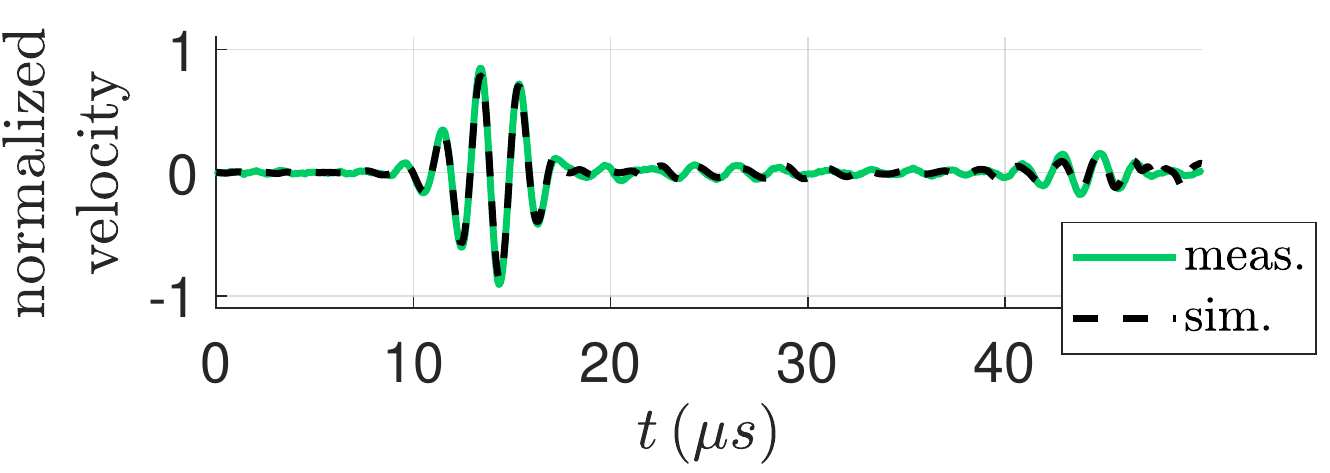}}
    \hfil
    \subcaptionbox{Plate with a notch of $\SI{0.8}{\milli\meter}$ depth.}
    {\includegraphics[width=\MyWidth]{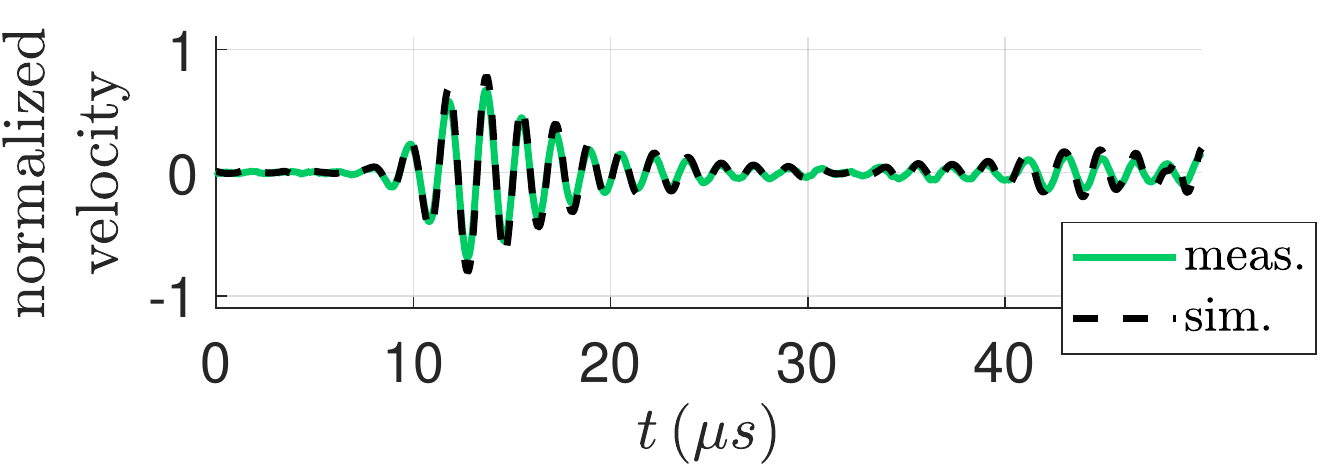}}
    \caption{
    Comparison of $z$-velocity at point $\vx=(\SI{-60}{\milli\meter},\SI{0}{\milli\meter})^\tr$ between forward models for reconstruction parameter $\vpara^{\min}$ (sim.) and the measurement (meas.).}
    \label{fig:Comparison}
\end{figure}

The optimization procedure should achieve the best possible agreement between simulation and measurement at the measurement points. To give an impression of the accuracy of the model, Figure~\ref{fig:Comparison} shows a direct comparison of the signals. Although the inverse method only uses the envelopes to compute the objective function, the phase of the signal from the simulation also matches the measured signal. The matching phase can be explained by the calculation of the transfer function, where the phase information is adjusted.

\begin{figure}[ht!]
    \centering
    \subcaptionbox{Plate with a notch of $\SI{0.4}{\milli\meter}$ depth.}
    {\includegraphics[width=\MyWidth]{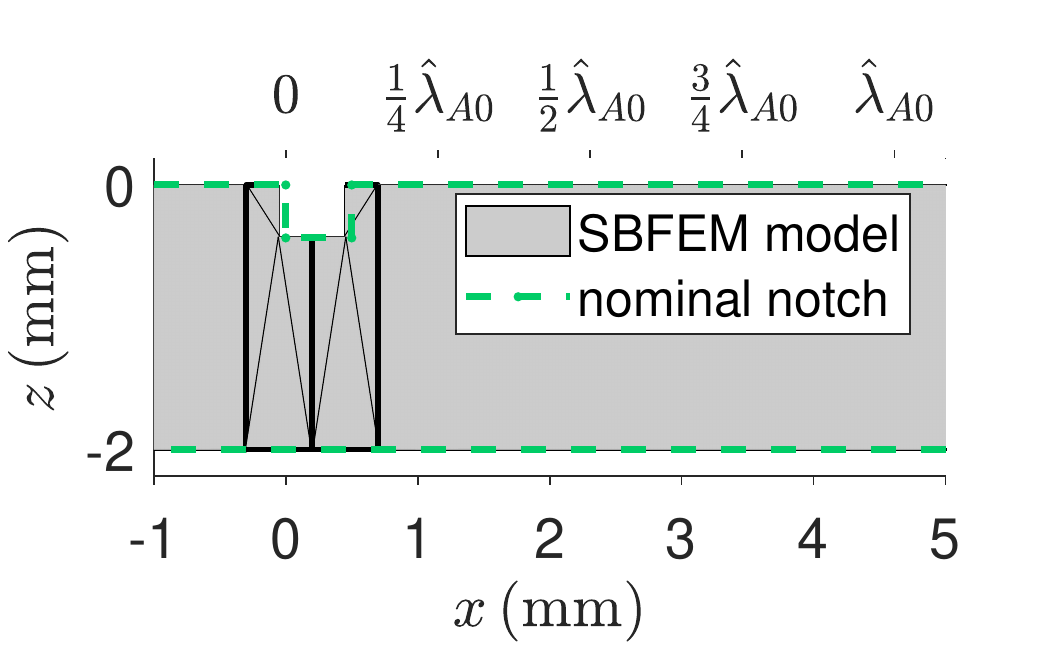}}
    \hfil
    \subcaptionbox{Plate with a notch of $\SI{0.8}{\milli\meter}$ depth.}
    {\includegraphics[width=\MyWidth]{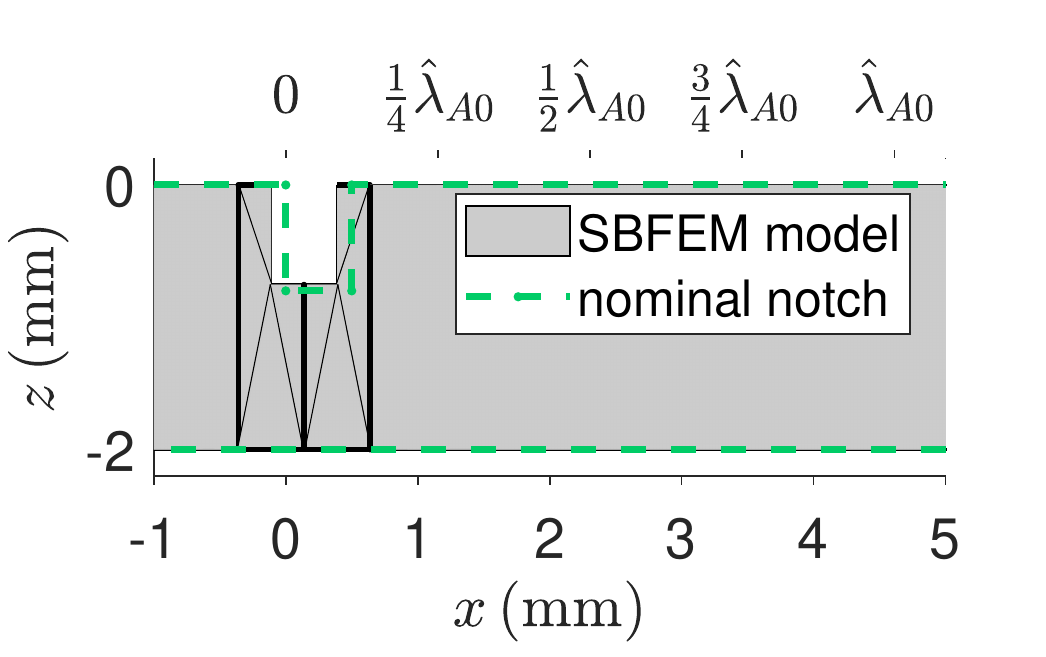}}
    \caption{Reconstruction model based on $\vpara^{\min}$ compared to the outline of the nominal notch based on $\vpara^{*}$, shown in green.}
    \label{fig:FinalReconstruction}
\end{figure}

For the evaluation of the reconstruction we refer to Figure~\ref{fig:FinalReconstruction}. In this figure the SBFEM mesh of the reconstruction is plotted against the nominal geometry shown in green.
The figures show only a small section around the notch. 
It should be noted that the actual parameter space is much larger.
The figure  shows a typically numbered $x$-axis at the bottom, while an $x$-axis measuring with the A0-wavelength of the center frequency $\ecf$ is shown at the top. This additional axis should make it possible to estimate the error for other plate dimensions, since the wavelength is the characteristic length of the problem.
The error for the x-position $|\para^*_1 - \para^{\min}_1|$ of the notch is $\SI{+0.05}{\milli\meter}$ and $\SI{+0.11}{\milli\meter}$ for the plate with a notch depth of $\SI{0.4}{\milli\meter}$ and $\SI{0.8}{\milli\meter}$, respectively.
In relation, the error for the notch depth $|\para^*_2 - \para^{\min}_2|$ is $\SI{+0.01}{\milli\meter}$ and $\SI{+0.05}{\milli\meter}$.
As already mentioned, the actual geometry may deviate from the nominal geometry due to manufacturing errors.

\FloatBarrier

\section{Conclusion}~\label{sec:Conclusion}

An algorithm has been presented that enables defect reconstruction using guided ultrasonic waves by solving an inverse optimization problem. The output of a simulation model, which contains a parametric representation of the defect, is intended to approximate the experimental data as closely as possible. For this purpose, a gradient-based method is used, where the computation of the derivatives according to the defect parameters is based on Algorithmic Differentiation. A second aspect of the presented work is the highly efficient Scaled Boundary Finite Element Method, which requires very small, final linear systems of equations to approximate the linear waveguide equation.

In this paper, which extends the theoretical results from~\cite{bulling2022defect} with experimental validation, the successful reconstruction of two notches in steel plates is presented. Although the reconstruction algorithm is formulated in a very general way, the initial focus was on the simplest case of wave propagation, where the ultrasonic wave can be assumed to be planar, and the problem can be reduced to a two-dimensional cross-sectional model. An error between nominal geometry and reconstruction in the order of 0.05 mm was achieved. Thus, it could be shown that, first, a geometric reconstruction of defects is possible and, second, that a cross-sectional model is capable of describing the wave field with sufficient accuracy for a reconstruction, given a suitable experimental setup.

The first important finding from this work is that the envelopes of the recorded data provide sufficient information while leading to an objective function that allows optimization with a gradient method. The second finding is that another linear fit must be introduced for the sensor's force distribution. This allows reconstruction even for sensors that are not fully characterized. Future work will extend the reconstruction algorithm to three-dimensional models of wave propagation. 

\section*{Acknowledgments}
The authors gratefully acknowledge the German Research Foundation for funding (DFG project number 428590437). In addition, the authors would like to thank Yevgeniya Lugovtsova and Tobias Homann for their support in the LDV measurement.

\printbibliography

\end{document}